\documentclass[12pt]{amsart}
\usepackage{amssymb,amscd,amsfonts,latexsym}   
\newtheorem{thm}{Theorem}[section] 
\newtheorem{prop}[thm]{Proposition}       
\newtheorem{cor}[thm]{Corollary}           
\newtheorem{lem}[thm]{Lemma}                   

\newtheorem{rem}[thm]{Remark}
\newtheorem{ex}[thm]{Example}

\numberwithin{equation}{section}

\def\R{{\bf R}} 
\def\Af{{\bf A}} 
\def\Pr{{\bf P}} 
\def\H{{\bf H}} 
\def\P{{\bf{\Psi}}}

\def\Ql{{\bf Q}_{\ell}}
\def\Zl{{\bf Z}_{\ell}}
\def\rnd{{\boldsymbol\cdot}}
\def\Q{{\bf Q}}
\def\C{{\bf C}}  
\def\Z{{\bf Z}}

\def\tsX{{\tilde X^*}}
\def\tsS{{\tilde S^*}}
\def\sX{{X^*}}

\newcommand{\bu}{{\bullet}}
\def\log{{\text{log}~}}

\newcommand{\ie}{{\it i.e.\/}\ }
\newcommand{\eg}{{\it e.g.\/}\ }
\newcommand{\cf}{{\it cf.\/}\ }
\newcommand{\resp}{{\it resp.\/}\ }
\newcommand{\op}{{\it op.cit.\/}\ }

\begin{document}

\title[The local monodromy as a generalized algebraic correspondence]{The local monodromy as a generalized\\\smallskip algebraic corresponcence}

\author[Caterina Consani]{Caterina Consani$^{\text{\dag}}$\\\medskip  with an appendix by Spencer Bloch}
\thanks{\noindent$^{\text{\dag}}$Partially supported by the NSF grant DMS-9701302}
\address{C.~Consani: Department of Mathematics \\
Massachusetts Institute of Technology \\
Cambridge, Ma.\ \ 02139 USA
}
\email{katia@math.mit.edu}
\address{S.~Bloch: Department of Mathematics \\
The University of Chicago \\
Chicago, Il.\ \ 606139 USA}

\email{bloch@math.uchicago.edu}
\maketitle

\section*{Introduction}

 Let $X$ be a proper and smooth variety over a local field $K$ and let $\mathcal X$ be a regular model of $X$ defined over the ring of integers $\mathcal O_K$ of $K$. When $\mathcal X$ is smooth over $\mathcal O_K$, the Tate conjecture equates the $\ell$--adic Chow groups of algebraic cycles on the geometric special fibre $X_{\bar k}$ of $\mathcal X \to \text{Spec}(\mathcal O_K)$ with the Galois invariants in $H^{2*}(X_{\bar K},\Ql(*))$. One of the results proved in \cite{C} (\cf Corollary 3.6) shows that the Tate conjecture for smooth and proper varieties over finite fields together with the monodromy--weight conjecture imply a generalization of the above result in the case of semistable reduction. Namely, let
$\wp \in \text{Spec}(\mathcal O_K)$ be a place over which the special fibre $\mathcal X \times \text{Spec}(k(\wp)) = Y$ is a reduced divisor with normal crossings in $\mathcal X$ (\ie semistable place). Then, assuming the above two conjectures, the $\ell$--adic groups of algebraic cycles modulo rational equivalence on the $r$--fold intersections of components of $Y$ ($r \ge 1$) are related with Galois invariant classes on the Tate twists $H^{2*-(r-1)}(X_{\bar K},\Ql(*-(r-1)))$.

An interesting case is when one replaces $X$ by $X \times_K X$, so that Galois invariant cycles may be identified with Galois equivariant maps $H^*(X_{\bar K},\Ql) \to H^*(X_{\bar K},\Ql(\cdot))$. Examples of such maps are the powers $N^i$ of the logarithm of the local monodromy  around $\wp$.  The operators $N^i: H^*(X_{\bar K},\Ql) \to H^*(X_{\bar K},\Ql(-i))$ determine classes $[N^i] \in H^{2d}((X \times X)_{\bar K},\Ql(d-i))$ ($d = \dim~X_{\bar K}$) invariant under the decomposition group. In this paper we study in detail the structure of $[N^i]$ when the special fibre $Y$ of $\mathcal X$ has at worst triple points as singularities. That is, we exhibit the corresponding algebraic cycles on the (normal crossings) special fibre $T = \cup_iT_i$ of a resolution $\mathcal Z$ of $\mathcal X \times_{\mathcal O_K} \mathcal X$.

 Denote by $\tilde N = 1 \otimes N + N \otimes 1$ the monodromy on the product, and let $F$ be the geometric Frobenius. Then the classes $[N^i]$ naturally determine elements in $\text{Ker}(\tilde N) \cap H^{2d}((X \times X)_{\bar K},\Ql(d-i))^{F=1}$. Assuming the monodromy--weight conjecture on the product (\ie the monodromy filtration $L_\rnd$ on $H^*((X \times X)_{\bar K},\Ql)$ coincides--up to a shift--with the filtration by the weights of the Frobenius \cf~\cite{RZ}), the following identifications hold  
\begin{equation}\label{eq0}
\text{Ker}(\tilde N) \cap H^{2d}((X \times X)_{\bar K},\Ql(d-i))^{F=1} \simeq \biggl ((gr^L_{2(d-i)}H^{2d}(T,\Ql))(d-i)\biggr )^{F=1} \simeq
\end{equation}
\[
\simeq \biggl [\frac{\text{Ker}(\rho^{(2(i+1))}: H^{2(d-i)}(\tilde{T}^{(2i+1)},\Ql)(d-i) \to H^{2(d-i)}(\tilde{T}^{(2(i+1))},\Ql)(d-i))}{\text{Image}~\rho^{(2i+1)}}\biggr ]^{F=1}.
\]\vspace{.1in}

Here $\tilde T^{(j)}$ denotes the normalization of the $j$--fold intersection on the closed fibre $T$.
These isomorphisms show that the classes $[N^i]$ have representatives in cohomology groups of some precise strata of $T$. Moreover, the Tate conjecture and the semisimplicity of the Frobenius for the smooth schemes $\tilde T^{(j)}$ would imply that these classes are algebraic. We refer to $\S~\ref{0}$, \eqref{diff} for the description of the restriction maps $\rho$ in \eqref{eq0}.

 To better understand the geometry related to the desingulatization process $\mathcal Z \to \mathcal X \times_{\mathcal O_K} \mathcal X$, and to avoid at first some technical complications connected to the theory of the nearby cycles in mixed characteristic, we start by investigating this
problem in equal characteristic zero (\ie for semistable degenerations over a disk). There, one can take full advantage of many geometric results based on the theory of the mixed Hodge structures. Under the assumption of the monodromy--weight conjecture and using some techniques of \cite{RZ}, our results generalize to mixed characteristic. The cycles we exhibit on $\tilde T^{(2i+1)}$ explain geometrically the presence of poles on specific local factors of the L--function related to the fibre product $X \times X$. In fact, theorem~\ref{th1} equates, under the assumption of the semisimplicity of the action of the Frobenius $F$ on the inertia invariants $H^{*}((X \times X)_{\bar K}, \Ql)^{I}$, the rank of any of the groups in \eqref{eq0} with $
\operatorname*{\text{ord}}_{s = d-i}\det(Id - FN(\wp)^{-s} | H^{2d}((X \times X)_{\bar K}, \Ql)^{I})$. Here, $N(\wp)$ denotes the number of elements of the finite field $k(\wp)$.

A study of the local geometry of the normal--crossings special fibre $T$ shows that $[N^i]$ are represented by certain natural ``diagonal cycles'' on $\tilde T^{(2i+1)}$ together with a cycle supported on the exceptional part of the stratum that arises because the classes $[N^i]$ must belong to the kernel of the restriction map $\rho^{(2(i+1))}$ (\cf \eqref{eq0}). This result is obtained via the introduction of a generalized correspondence diagram for the map
\begin{equation}\label{eq1}
N^i: \H^*(Y, gr^L_{r+i}\R\P(\Q_{\mathcal X})) \to  \H^*(Y, (gr^L_{r-i}\R\P(\Q_{\mathcal X}))(-i)).
\end{equation}

This morphism describes the monodromy action on the $E_1$--term of the spectral sequence of weights for the filtered complex of the nearby cycles  $(\R\P(\Q_{\mathcal X}),L_\rnd)$ (\cf $\S~\ref{a}$, \eqref{2}). For $i > 0$, the classes $[N^i]$ do not describe an algebraic correspondence in the classical sense. In fact, the algebraic cycles representing them are only supported on higher strata of the special fibre $T$ (\ie on $\tilde T^{(2i+1)}$) and they do not naturally determine classes in the cohomology of $T$. This is a consequence of the fact that for $i > 0$, the cocycle $[N^i]$ does not have weight zero in the $\ell$--adic cohomology of the fibre product $(X \times X)_{\bar K}$, as one can easily check from \eqref{eq0}. Nonetheless, we expect that each of these classes  supplies a refined information on the degeneration. Namely, we conjecture that the geometric description that we obtain up to triple points can be generalized to any kind of semistable singularity via a thorough combinatoric study of the toric singularities of the special fibre of the fibre product resolution $\mathcal Z$.

The correspondence diagram related to the map \eqref{eq1} is 
 built up from the hypercohomology of the Steenbrink filtered resolution $(A^\bu_{\mathcal X},L_\rnd)$ of 
  $\R\P(\Q_{\mathcal X})$. In $\S~\ref{b}$ we establish the necessary  functoriality properties of 
the Steenbrink complex and its $L_\rnd$--filtration. A difficult point in the description of the correspondence diagram is related to the definition of a product structure on the $E_1$--terms of the spectral sequence of weights. In fact, in the Example~\ref{ex1} we show that there is no canonical definition of a product structure for $(A^\bu_{\mathcal X},L_\rnd)$ in the filtered category. Equivalently said, the monodromy filtration $L_\rnd$ is not multiplicative on the level of the filtered complexes. A partial product, canonical only on the $E_2 = E_{\infty}$-terms is provided in  the Appendix. This suffices for purposes of our paper. 

\subsection*{\it Acknowledgments} I would like to thank Kazuya Kato for suggesting the study of the local monodromy as a Galois invariant class, eventually algebraic. I also would like to acknowledge an interesting conversation with Alexander Beilinson on some of the arguments presented here. I gratefully thank the Institut des Hautes \'Etudes Scientifiques for the kind hospitality received during my stay there in January 1997. Finally, it is my pleasure to thank Spencer Bloch for his constant support and for many fruitful suggestions  I received from him during the preparation of this paper.
\vspace{.2in}

\section{Notations and techniques from mixed Hodge theory}\label{0}

In this paragraph we introduce the main notations and recall some results on the mixed Hodge theory of a degeneration.\vspace{.1in}

We denote by $X$ a connected, smooth, complex analytic manifold and we let $S$ be the unit disk. We write $f: X \to S$ for a proper, surjective morphism and we let $Y = f^{-1}(0)$ be its special fibre. We assume that $f$ is {\it smooth} at every point of $\sX = X \smallsetminus Y$ and that the special fibre $Y$ is an algebraic divisor with {\it normal--crossings}. The local description of $f$ near a closed point $y \in Y$ is given by:
\[
f(z_1,\ldots,z_m) = z_1^{e_1}\cdots z_k^{e_k}
\]

\noindent for $k \le m = \dim X$ and $\{z_1,\ldots,z_m\}$ a local coordinate system on a neighborhood of $y$ in $X$ centered at $y$ and $e_i \in \Z,~e_i \ge 1$. The fibres of $f$ have then dimension $d = m - 1$.

We {\it fix a parameter} $t \in S$. For $t \neq 0$, let $f^{-1}(t) = X_t$ be the fibre at $t$. Because the restriction of $f$ at $S^* = S \smallsetminus \{0\}$ is a $C^{\infty}$, locally trivial fibre bundle, the positive generator of $\pi_1(S^*,t) \simeq \Z$ induces an automorphism $T_t$ of $H^*(X_t,\Z)$, called the local monodromy. We will always assume throughout the paper that $T_t$ is {\it unipotent} (\ie $(T_t-1)^{i+1} = 0$, on $H^i(X_t,\Z)$). This condition is naturally verified when for example  $\text{g.c.d.}(e_i,~i\in[1,k]) = 1$, $\forall y \in Y$. Under these assumptions, the logarithm of the local monodromy is defined to be the finite sum:
\[
N_t := \log T_t = (T_t - 1) - \frac{1}{2}(T_t - 1)^2 + \frac{1}{3}(T_t - 1)^3 - \cdots
\]

It is known (\cf~\cite{D2}) that the automorphisms $T_t$ of $H^i(X_t,\C)$ ($t \in S^*$), are the fibres of an automorphism $T$ of the fibre bundle $\R^if_*(\Omega^\bu_{X/S}(\log Y))$ over $S$, whose fibre at $0$ is described as $T_0 = \exp(-2\pi iN_0)$.
By definition, the endomorphism $N_0$ is  the residue at $0$ of the Gauss-Manin connection $\nabla$ on the ``canonical prolongation'' $\R^if_*(\Omega^\bu_{X/S}(\log Y))$ of the locally free sheaf $\R^if_*(\Omega^\bu_{X^*/S^*})$. 
Because of the definition of $T_0$, it makes sense to think of a nilpotent map $N := - \frac{1}{2\pi i}~\log T$ as the monodromy operator on the degeneration $f: X \to S$. Via the canonical isomorphism (\cf~\cite{S}, Thm. 2.18)($t \in S$): 
\begin{equation*}
\R^if_*(\Omega^\bu_{X/S}(\log Y)) \otimes_{\mathcal O_{S}} k(t) \overset{\simeq}\to \H^i(X_t,\Omega^\bu_{X/S}(\log Y) \otimes_{\mathcal O_X} \mathcal O_{X_t})
\end{equation*}\noindent
where $k(t)$ is the residue field of $\mathcal O_{S}$ at $t$, we can see the map $N_0$ as an endomorphism of the hypercohomology of the relative de Rham complex $\Omega^\bu_{X/S}(\log Y) \otimes_{\mathcal O_X} \mathcal O_Y$. This complex represents in the derived category  $D^+(Y,\C)$ of the abelian category of sheaves of $\C$--vector spaces on $Y$, the complex of the nearby cycles $\R\P(\C)$. Namely, there exists a non-canonical quasi-isomorphism (\ie depending on the choice of the parameter $t$ on $S$) $\Omega^\bu_{X/S}(\log Y) \otimes_{\mathcal O_X} \mathcal O_Y \simeq \R\P(\C_{\tsX}) := i^{-1}\R k_*\C_{\tsX}$ (\cf~\cite{S}, $\S~2$). We refer to the following commutative diagram for the description of the maps: 
\begin{equation*}
\begin{CD}
\tsX @>k>> X @<i<< Y\\
@VVV @VV{f}V @VVV\\
\tsS @>p>> S @<<<\{0\}.
\end{CD}
\end{equation*} 

The space $\tsS = \{u \in \C~| \text{Im}~u > 0\}$ is the upper half plane, the map $p: \tsS \to S$ $p(u) = \exp(2\pi i u) = t$, makes $\tsS$ in a universal covering of $S^*$ and $\tsX$ is the pullback $X \times_S \tsS$ of $X$ along $p$. The morphism $k$ is the natural projection. It factorizes through $\sX$ by means of the injection $j: \sX \to X$. Finally, $i$ is the closed embedding.

Steenbrink defined  a mixed Hodge structure on the hypercohomology of the unipotent factor of the complex of the nearby cycles  $\H^*(X,\Omega^\bu_{X/S}(\log Y) \otimes_{\mathcal O_X} \mathcal O_{Y^{\text{red}}})$. This is frequently referred as the  limiting mixed Hodge structure. 

We will assume from now on that $f$ is {\it projective}. Then, the weight filtration on the limiting mixed Hodge structure is the one induced by the nilpotent endomorphism $N$, namely by the logarithm of the unipotent Picard-Lefschetz transformation $T$ that is already defined at the $\Q$-level. This filtration, which one usually refers to as the monodromy--weight filtration $L_\rnd$, is defined inductively. On the limiting cohomology $H^i(\tsX,\Q)$, it is increasing and has lenght at most $2i$. By the local monodromy theorem $N^{i+1} = 0$, hence one sets $L_0 = \text{Im}~N^i$ and $L_{2i-1} = \text{Ker}~N^i$. The monodromy filtration $L_\rnd$ becomes  a convolution product of the kernel and the image filtration relative to the endomorphism $N$. These filtrations are defined as
\[
K_l~H^i(\tsX,\Q) := \text{Ker}~N^{l+1},\qquad I^j~H^i(\tsX,\Q) := \text{Im}~N^j
\]

\noindent and their convolution is
\begin{equation}\label{monodromy}
L = K * I, \qquad L_k := \sum_{l-j = k}K_l \cap I^j.
\end{equation}

It is a very interesting fact that there is no explicit construction of the monodromy-weight filtration $L_\rnd$ on $\Omega^\bu_{X/S}(\log Y) \otimes_{\mathcal O_X} \mathcal O_Y$ itself. The  filtration $L_\rnd$ is  defined on a complex $A^\bu_{\C}$ which is a resolution of $\Omega^\bu_{X/S}(\log Y) \otimes_{\mathcal O_X} \mathcal O_{Y^{\text{red}}}$. More precisely, the complex $\Omega^\bu_{X/S}(\log Y) \otimes_{\mathcal O_X} \mathcal O_{Y^{\text{red}}}$ is isomorphic, in the derived category $D^+(Y,\C)$, to the complex $A^\bu_{\C}$ of $\mathcal O_X$--modules supported on $Y$. The complex $A^\bu_{\C}$ is the simple complex associated to the double complex ($p,~q \ge 0$):
\[
A^{p,q}_{\C} := \Omega_X^{p+q+1}(\log Y)/W_q\Omega_X^{p+q+1}(\log Y) 
\]

\noindent where $W_*\Omega^\bu_X(\log Y)$ is the weight filtration by the order of log-poles (\cf~\cite{D},~$\S~3$). The differentials on it are defined as follows
\[
d': A^{p,q}_{\C} \to A^{p+1,q}_{\C},\qquad d'(\omega) = d\omega
\]

\noindent is induced by the differentiation on the complex $\Omega^\bu_X(\log Y)$ and 
\[
d'': A^{p,q}_{\C} \to A^{p,q+1}_{\C},\qquad d''(\omega) = \theta \wedge \omega
\]

\noindent where $\theta := f^*(\frac{dt}{t}) = \sum_{i = 1}^ke_i\frac{dz_i}{z_i}$ is the form definying the quasi-isomorphism we mentioned before (\cf~\cite{S}, $\S~4$)
\begin{equation*}
\Omega^\bu_{X/S}(\log Y) \otimes_{\mathcal O_X} \mathcal O_Y \overset{\wedge~\theta}\to A^\bu_{\C}.
\end{equation*}

 The total differential on $A^\bu_{\C}$ is $d = d' + d''$. The weight filtration $W_*\Omega_X^\bu(\log Y)$ induces a corresponding filtration on $A^\bu_{\C}$ ($r \in \Z$):
\begin{equation}\label{weight}
\begin{CD}
W_rA^{p,q}_{X,\C} =:  W_{r+q+1}\Omega^{p+q+1}_X(\log Y)/W_q\Omega^{p+q+1}_X(\log Y).
\end{CD}
\end{equation}

The filtration that $W_rA^\bu_{\C}$ induces on $\H^*(Y,A^\bu_{\C}) \simeq \H^*(\tsX,\C)$ is the kernel filtration $K$ (\cf \eqref{monodromy})
\[
K_rH^*(\tsX,\C) = W_r\H^*(Y,A^\bu_{\C}) =: \text{Im}\biggl (\H^*(Y,W_rA^\bu_{\C}) \to \H^*(Y,A^\bu_{\C})\biggr ) = \text{Ker}~N^{r+1}.
\]

The monodromy-weight filtration is then defined as
\[
L_rA^{p,q} := W_{2q+r+1}\Omega^{p+q+1}_X(\log Y)/W_q\Omega^{p+q+1}_X(\log Y).
\]

Via Poincar\'e residues, the related graded pieces have the following description
\begin{equation}\label{1}
gr^L_rA^\bu_{\C} \simeq \bigoplus_{k \ge \text{max}(0,-r)}(a_{2k+r+1})_*\Omega_{\tilde Y^{(2k+r+1)}}^\bu[-r-2k].
\end{equation}

Here, we have denoted by  $\tilde Y^{(m)}$ the disjoint union of all intersections $Y_{i_1} \cap \ldots \cap Y_{i_m}$ for $1 \le i_1 < \ldots < i_m \le n$ ($Y = Y_1 \cup \ldots \cup Y_n$). We write $(a_m)_*: \tilde Y^{(m)} \to X$ for the natural projection.

The monodromy operator $N$ is induced by an endomorphism $\tilde\nu$ of $A^\bu_{\C}$ which is defined as $(-1)^{p+q+1}$ times the natural projection
\[
\nu: A^{p,q}_{\C} \to A^{p-1,q+1}_{\C}.
\]

The endomorphism $\tilde\nu$ is characterized by its behavior on the $L$-filtration, namely 
\[
\tilde\nu(L_rA^\bu_{\C}) \subset L_{r-2}A^\bu_{\C}
\]
and the induced map
\begin{equation}\label{nu}
\tilde\nu^r: gr^L_rA^\bu_{\C} \to gr^L_{-r}A^\bu_{\C}
\end{equation}

\noindent is an isomorphism for all $r \ge 0$. The complex $A^\bu_{\C}$ contains the subcomplex $W_0A^\bu_{\C} = \text{Ker}(\tilde\nu)$ that is known to be a resolution of $\C_{Y}$. The filtration $L$ and the Hodge filtration $F$ on $A^\bu_{\C}$ induce \resp the kernel and $F$ filtration on $W_0A^\bu_{\C}$. The resulting mixed Hodge structure on $H^*(Y,\C)$ is the canonical one. Similarly, the homology $H_*(Y,\C)$ (\ie $H^*_{Y}(X,\C)$) with its mixed Hodge structure is calculated by the hypercohomology of the complex $\text{Coker}(\tilde\nu)$.

Because of the description given in \eqref{1}, the spectral sequence of hypercohomology of the filtered complex $(A^\bu_{\C},L)$ (frequently referred as the {\it weight spectral sequence} of $\R\P(\C)$) has the $E_1$ term  explicitly given by
\begin{equation}\label{1a}
\begin{split}
E_1^{-r,n+r} &= \bigoplus_{k \ge \text{max}(0,-r)}H^{n-r-2k}(\tilde Y^{(2k+r+1)},\C)\\
d_1 &= \sum_k ((-1)^{r+k}d_1' + (-1)^{k-r}d_1'')
\end{split}
\end{equation}

\noindent where
\begin{equation}\label{diff}
\begin{split}
d_1' &= \rho^{(r+2k+2)} = \sum_{u = 1}^{r+2k+2}(-1)^{u-1} \rho_u^{(r+2k+2)}\\
d_1'' &= -\gamma^{(r+2k+1)} = \sum_{u = 1}^{r+2k+1}(-1)^u \gamma_u^{(r+2k+1)}
\end{split}
\end{equation}

\noindent are the differentials and 
\begin{align*}
\rho_u^{(r+2k+2)} &= (\delta_u^{(r+2k+2)})^*: H^{n-r-2k}(\tilde Y^{(2k+r+1)},\C) \to H^{n-r-2k}(\tilde Y^{(2k+r+2)},\C)\\
\gamma_u^{(r+2k+1)} &= (\delta_u^{(r+2k+1)})_!: H^{n-r-2k}(\tilde Y^{(2k+r+1)},\C) \to H^{n-r-2k+2}(\tilde Y^{(2k+r)},\C)
\end{align*}

are \resp the restrictions and the Gysin maps induced by the inclusions ($u,t\in \Z$)
\[
\delta_u^{(t)}: Y_{i_1} \cap\cdots\cap Y_{i_t} \to Y_{i_1} \cap\cdots\cap (Y_{i_u})^{\hat{}}\cap\cdots\cap Y_{i_t}.
\]

Notice that the weight spectral sequence \eqref{1a} is built up from a filtered {\it double} complex, so that its $E_1$--term becomes a total complex of a double complex. This property distinguishes this weight spectral sequence from others as \eg the spectral sequence of weights which defines the mixed Hodge structure on a quasi--projective smooth complex variety (\cf \cite{D}). 

The complex $A^\bu_{\C}$ is the complex part of a cohomological mixed Hodge complex $A^\bu_{\Q}$ whose definition is less explicit than $A^\bu_{\C}$ and for which we refer to \cite{I}. This rational complex induces on $H^\rnd(\tsX,\Q)$ a rational mixed Hodge structure. The rational representative of the above spectral sequence \eqref{1a} is
\begin{equation}\label{1b} 
E_1^{-r,n+r} = \bigoplus_{k \ge \text{max}(0,-r)}H^{n-r-2k}(\tilde Y^{(2k+r+1)},\Q)(-r-k).
\end{equation}

The index in the round brackets outside the cohomology refers to the Tate twist. 
Both these spectral sequences degenerate at $E_2 = E_{\infty}$ and  they converge to $H^n(\tsX,\C)$ and $H^n(\tsX,\Q)$ respectively.

For curves (\ie $d = 1$), the degeneration of the weight spectral sequence provides the exact sequences
\begin{equation*}
0 \to E_2^{-1,2} \to H^0(\tilde Y^{(2)},\Q)(-1) \overset{d_1^{-1,2}} \to 
 H^2(\tilde Y^{(1)},\Q) \to H^2(\tsX,\Q) \to 0
\end{equation*}

\noindent and

\begin{equation}\label{1c}
0 \to H^0(\tsX,\Q) \to H^0(\tilde Y^{(1)},\Q) \overset{d_1^{0,0}} \to H^0(\tilde Y^{(2)},\Q) \overset{\alpha} \to H^1(\tsX, \Q).
\end{equation}

The differentials $d_1^{-1,2}$ and $d_1^{0,0}$ are defined as in \eqref{diff} and the map $\alpha$ in \eqref{1c} is the edge map in the spectral sequence. We also have a non canonical decomposition
$$
H^1(\tsX,\Q) = H^1(\tilde Y^{(1)},\Q) \oplus E_2^{-1,2} \oplus E_2^{1,0}.
$$

\noindent with $E_2^{1,0} = \text{Im}(\alpha)$.\vspace{.1in}

Steenbrink proves that the filtration induced on the abutment of the spectral sequence of the nearby cycles is the Picard-Lefschetz filtration, hence it is uniquely described by the following properties
\[
N(gr^L_{n+r}H^n(\tsX,\Q)) \subset (gr^L_{n+r-2}H^n(\tsX,\Q))(-1)
\]
and 
\begin{equation*}
N^r: gr^L_{n+r}H^n(\tsX,\Q) \overset{\simeq}\to (gr^L_{n-r}H^n(\tsX,\Q))(-r)
\end{equation*}

\noindent for $r > 0$. In the rest of the paper we will refer to it as the {\it monodromy filtration}.\vspace{.2in}

\section{The monodromy operator as algebraic cocycle}\label{a}

We keep the notations introduced in the last paragraph. As $n$ varies in $[0,2d]$ ($d =$ dimension of the fibre of $f: X \to S$) and $i \ge 0$, the power maps
\[
N^i: H^n(\tsX,\Q) \to H^n(\tsX,\Q)(-i)
\]

\noindent induced by the endomorphism $N: \R^nf_*(\Omega^\bu_{X/S}(\log Y)) \to \R^nf_*(\Omega^\bu_{X/S}(\log Y))$, define elements
$$
N^i \in \text{Hom}(H^\rnd(\tsX,\Q), H^\rnd(\tsX,\Q)(-i)) 
$$

\noindent which are invariant for the action of the local monodromy group $\pi_1$. They can be naturally identified with
\[
N^i \in \bigoplus_{n \ge 0}\biggr [H^{2d-n}(\tsX,\Q)(d) \otimes H^n(\tsX,\Q)(-i)\biggr ]^{\pi_1} = \biggr [H^{2d}(\tsX \times \tsX,\Q)(d-i)\biggl  ]^{\pi_1}.
\]

The space $\tsX \times_S \tsX$ is the generic fibre of the product degeneration $X \times_S X \to S$. After a suitable sequence of blow-ups  along $\text{Sing}(Y \times Y) \supset \text{Sing}(X \times_S X)$: 
\[
Z \to \cdots \to X \times_S X \to S
\]

 \noindent we obtain a normal--crossings degeneration $h: Z \to S$ with $Z$ non singular and whose generic fibre  is still $\tsX \times \tsX$. Its special fibre $T = {h}^{-1}(0) = T_1 \cup \cdots \cup T_N$ has normal crossings singularities. The local description of $h$ along $T$ looks like:
\[
h(w_1,\ldots,w_{2m}) = w_1^{e_1}\cdots w_r^{e_r}
\]
for $\{w_1,\ldots,w_{2m}\}$ a set of local parameters on $Z$ and $e_1,\ldots,e_r$ non-negative integers.\vspace{.1in} 

The semistable reduction theorem assures that modulo extensions of the basis $S$ and up to a suitable sequence of blow-ups and down along subvarieties of the special fibre $T$, we may eventually obtain from $h$ a semistable degeneration $W \to S$ with $W_0 = W_{0_1} \cup \ldots \cup W_{0_M}$ as special fibre. \vspace{.1in}

Because of the assumption of the unipotency of the local monodromy on $H^*(X_t,\C)$ (\cf~$\S$~\ref{0}), the local monodromy $\sigma$ of $h$ will be also unipotent. We then call $\tilde N = \log(\sigma)$. By the K\"unneth decomposition it results: $\tilde N = 1 \otimes N + N \otimes 1$ and we have:
\[
N^i \in \biggl (H^{2d}(\tsX \times \tsX,\Q(d-i))\biggl )^{\pi_1} = \text{Ker}(\tilde N) \cap H^{2d}(\tsX \times \tsX,\Q(d-i)).
\]

Let consider the monodromy filtration $L_\rnd$ relative to the degeneration $h$. We denote by $Hom_{MH}(\Q(0),V)$ ($Hom(\Q,V)$ shortly) the subgroup of  Hodge cycles of pure weight $(0,0)$ of a bifiltered $\Q$--vector space $V$: $(V,L,F)$, endowed with the corresponding mixed Hodge structure. Then, we have the following

\begin{prop}\label{prop1} For $i \ge 1$

 $N^i \in Hom_{MH}\biggl (\Q(0),\text{Ker}(\tilde N) \cap H^{2d}(\tsX \times \tsX,\Q(d-i))\biggr ) \simeq$
\[
\simeq Hom_{MH}\biggl (\Q(0),(gr^L_{2(d-i)}H^{2d}(\tsX \times \tsX,\Q))(d-i)\biggr ) \simeq 
\]
\[ 
\simeq Hom_{MH}\biggl (\Q(0),(gr^L_{2(d-i)}H^{2d}(T,\Q))(d-i)\biggr ) \simeq 
\]
\begin{equation*}
\simeq Hom\biggl (\Q,\frac{\text{Ker}(\rho^{(2(i+1))}: H^{2(d-i)}(\tilde{T}^{(2i+1)},\Q)(d-i) \to H^{2(d-i)}(\tilde{T}^{(2(i+1))},\Q)(d-i))}{\text{Image}~\rho^{(2i+1)}}\biggr ).
\end{equation*}

Here $\rho$ is the restriction map on cohomology and by $\tilde{T}^{(j)}$ we mean the disjoint union of all ordered $j$--fold intersections of the components of $T$ (\cf~$\S\ref{0}$).
\end{prop}

\begin{proof} The identification of $N^i$ with a Hodge cycle is a consequence of $N$ being a morphism in the category of Hodge structures. The first isomorphism derives from the well known facts that $\text{Ker}(\tilde N)$ has monodromic weight zero and that the Hodge cycles belong only to the graded piece $(gr^L_{2(d-i)}H^{2d}(\tsX \times \tsX,\Q))(d-i)$ of $\text{Ker}(\tilde N) \cap \bigoplus_j (gr^L_j H^{2d}(\tsX \times \tsX,\Q))(d-i)$. The second isomorphism comes from the local invariant cycle theorem, namely from the following exact sequence of pure Hodge structures (\cf~\cite{C}, lemma~3.3 and corollary~3.4)
\begin{equation*}
0 \to gr^L_{2(d-i)}H^{2d}(T,\Q) \to gr^L_{2(d-i)}H^{2d}(\tsX \times \tsX,\Q) \overset{N} \twoheadrightarrow gr^L_{2(d-i-1)}H^{2d}(\tsX \times \tsX,\Q)(-1)
\end{equation*}

 Finally, the last isomorphism is a consequence of the description of the graded piece $(gr^L_{2(d-i)}H^{2d}(T,\Q))(d-i)$ as sub--Hodge structure of $(gr^L_{2(d-i)}H^{2d}(\tsX \times \tsX,\Q))(d-i)$ (\cf\op~lemma~3.3).\end{proof}\vspace{.1in}

Proposition~\ref{prop1} shows how the operators $N^i$ can be detected by  classes $[N^i]$ in the cohomology of a fixed stratum of the special fibre $T$. Equivalently, we can say that $N^i$ determine classes $[N^i] \in \H^{2d}(T,(gr^L_{-2i}\R\P_{h}(\Q))(d-i))$ in the $(E_1^{2i,2(d-i)})(d-i)$-term of the spectral sequence of weights for the degeneration $h$. Here we write $gr^L_{-2i}\R\P_{h}(\Q)$ for $gr^L_{-2i}A^\bu_{W,\Q}$.\vspace{.1in} 

The goal of this paper is to identify the class $[N^i]$ with an algebraic cocycle related to the degeneration $f: X \to S$. In all those cases that we will consider in the paper, this identification is obtained via  a  ``correspondence-type'' map ($i \ge 0$)
\[
N^i: \H^*(Y, gr^L_rA^\bu_{X,\Q}) \to \H^*(Y, (gr^L_{r-2i}A^\bu_{X,\Q})(-i)) = \H^*(Y, gr^L_r(A^\bu_{X,\Q}(-i)))
\]

\noindent which makes the following diagram  commute
\begin{equation}\label{2}
\begin{CD}
\H^*(T,gr^L_rA^\bu_{Z,\Q}) @>{[N^i]~\cdot}>> \H^{2d+*}(T, (gr^L_{r-2i}A^\bu_{Z,\Q})(d-i)) = E_1^{-r+2i,*+r+2(d-i)}\\
@A{(p_1)^*}AA @VV{(p_2)_*}V \\
E_1^{-r,*+r} = \H^*(Y, gr^L_rA^\bu_{X,\Q}) @>{N^i}>>  \H^*(Y, gr^L_{r}(A^\bu_{X,\Q}(-i))) = E_1^{-r+2i,*+r-2i}
\end{CD}
\end{equation}\vspace{.1in}

The projections $p_1,~p_2: \tsX \times \tsX \to \tsX$ on the first and second factor, determine pullbacks and pushforwards on the hypercohomology as we shall describe in $\S$~\ref{b}.

From the theory we will explain in the next paragraphs and in the Appendix it will follow that $N^i$ has the expected shape. Namely, it is zero when $N^i = 0$ and it is the identity when $N^i$ induces an isomorphism on $E_2^{-r,*+r}$. Also, it will result that $p_1^*$, $(p_2)_*$ and $[N^i]~\cdot$ all commute with the differential on $E_1$. That will imply an induced commutative diagram on $E_2$. 

For $i = 0$, \ie when the correspondence map is the identity, proposition~\ref{prop1} can be slightly generalized, using the theory developed in \cite{C} (\cf~lemma~3.3 and corollary~3.4) and in \cite{BGS} so that the identity operator is seen as an element in 
\begin{multline*}
Hom_{MH}\biggl (\Q,\frac{\text{Ker}(\rho^{(2)}: H^{2d}(\tilde{T}^{(1)},\Q)(d) \to H^{2d}(\tilde{T}^{(2)},\Q)(d))}{\text{Im}~(-i^*\cdot i_*: H_{2(d-1)}(T^{(1)},\Q)(d-1) \to H^{2d}(\tilde{T}^{(1)},\Q)(d)}\biggr ) \simeq\\\\ \simeq Hom_{MH}\biggl (\Q,\text{Im}(i^*: H^{2d}(T,\Q)(d) \to H^{2d}(\tilde{T}^{(1)},\Q)(d))\biggl ).
\end{multline*}

Here the map $i^*$ (\resp $i_*$) represents the pullback (\resp pushforward)  relative to the embedding $T^{(1)} \to T$. Proposition~\ref{prop1} shows this class as a Hodge cocycle in $H^{2d}(\tsX \times \tsX,\Q(d))$. That agrees with the classical theory of algebraic correspondences describing the identity map via an algebraic correspondence with the cycle diagonal. Namely, the identity is determined by the diagonal $\Delta_{\tsX} \subset \tsX \times \tsX$ seen as specialization of the cycle diagonal on $\mathcal X \times \mathcal X$ on the fibre product $\tsX \times \tsX$. (\cf~\cite{K}).  

The cases described in the next paragraphs will also supply some evidence for our expectation that $[N^i]$ can be always described by an algebraic (motivic) cocycle. Finally, notice that the calculation on the $E_1$ involves the cohomology of individual components of the strata and it is therefore in some sense local, whereas $E_2$ introduces relations among components of strata, so that any calculation on it becomes of global nature. That is the reason why the description of the monodromy cycle is carried out mainly at a local level in this paper. \vspace{.2in}

\section{Functoriality of the Steenbrink complex and remarks on products}\label{b}

Let $g: Z \to X$ be a morphism between two connected, complex analytic manifolds over a disk $S$. Let $f: X \to S$ and $h: Z \to S$ be the degeneration maps. Let assume that both $Z$ and $X$ are smooth over $\C$ and  they have algebraic special fibres $f^{-1}(0) = Y$ and $h^{-1}(0) = T$ with normal crossings. We have the following commutative diagram
\begin{equation*}
\begin{array}{ccccc}
T &&{\longrightarrow} && Y\\
\downarrow{i'} &&&& \downarrow{i} \\ 
Z &&\stackrel{g}{\longrightarrow} && X \\
&{h}\searrow && \swarrow{f} &\\ && S &&.
\end{array}
\end{equation*}

Locally on the special fibres, $f$ and $h$ have the following description
\[
f(z_1,\ldots,z_m) = z_1^{e_1}\cdots z_k^{e_k};\qquad h(w_1,\ldots,w_M) = w_1^{e'_1}\cdots w_K^{e'_K}
\]

\noindent for $\{z_1,\ldots,z_m\}$ and $\{w_1,\ldots,w_M\}$ local parameters \resp on $X$ and $Z$, $1 \le k \le m,~1 \le K \le M$ and $e_1,\ldots,e_k;~e'_1,\dots,e'_K$ integers.\vspace{.1in}

Because $g^{-1}(Y) = T$, at any point $y \in g(T) \subset Y$ ($y = g(t)$, for some $t \in T$) where the local description of $Y$ is $z_1^{e_1}\cdots z_k^{e_k} = 0$, the pullback sections $g^*(z_{i_j})$ ($\forall~1 \le i_j \le k$) define divisors on $Z$ supported on $T$ (not necessarily reduced or irreducible). 

Let order the components of $Y$ as $Y = Y_1 \cup\ldots \cup Y_k$ and let denote by $\tilde Y^{(r)}$ the disjoint union of all intersections $Y_{i_1}\cap\ldots \cap Y_{i_r}$ for $1 \le i_1 <\cdots < i_r \le k$. There is a local system $\epsilon$ of rank one on $\tilde Y^{(r)}$ of standard orientations of $r$ elements (\cf~\cite{D}). The canonical morphism
\[
g^*\Omega^\bu_X(\log Y) \to \Omega^\bu_Z(\log T)
\]

\noindent is a map of bifiltered complexes with respect to the weight and the Hodge filtrations on $X$ and $Z$ \resp (\cf~\op). In particular it induces the following map of bicomplexes of sheaves supported on the special fibres ($r \ge 0$)
\[
g^*(W_rA^\bu_{X,\C}) \to W_rA^\bu_{Z,\C}
\]

\noindent where $A^\bu_{\C}$ is the Steenbrink complex which represents in the derived category the maximal subobject of the complex of nearby cycles where the action of the monodromy is unipotent (\cf~$\S~\ref{0}$). $W_rA^\bu_{\C}$ is the induced weight filtration on $A^\bu_{\C}$ (\cf~\eqref{weight}). Because the weight filtration on the complex $A^\bu_{\C}$ is induced by the weight filtration on the de Rham complex with log-poles, $g$ induces a map in the derived category
\[
g^*(W_r\R\P_f(\Q_X)) \to W_r\R\P_{h}(\Q_Z).
\]

Notice that $g^*(\frac{dz_{i_j}}{z_{i_j}}) \in W_1\Omega_Z^1(\log T)$, \ie pullbacks preserve poles. Hence, we deduce the functoriality of the monodromy filtration 
\[
g^*(L_rA^\bu_{X,\C}) \to L_rA^\bu_{Z,\C}.
\]
 
Because $g^{-1}$ is an exact functor, $g$ determines on the graded pieces a pullback map
\[
g^*: gr^L_rA^\bu_{X,\C} \to gr^L_rA^\bu_{Z,\C}
\]

\noindent where 
\[
gr^L_rA^\bu_{Z,\C} \simeq \bigoplus_{k \ge \text{max}(0,-r)}(a_{2k+r+1})_*\Omega_{\tilde T^{(2k+r+1)}}^\bu(\epsilon^{2k+r+1})[-r-2k].
\]

The functor $g^{-1}$ is also compatible with both differentials $d'$ and $d''$ on $A^\bu_{\C}$. Hence, $g^*$ induces a morphism of bifiltered mixed Hodge complexes ($F^\rnd =$ Hodge filtration \cf~\cite{D})
\[
g^*: (A^\bu_{X,\C},L,F) \to (A^\bu_{Z,\C},L,F)
\]

\noindent which in turn induces a map between the spectral sequences of weights
\[
g^*: E_1^{-r,q+r}(X) = \H^q(Y, gr^L_rA^\bu_X) \to \H^q(T, gr^L_rA^\bu_Z) = E_1^{-r,q+r}(Z).
\]

On the rational level this morphism between spectral sequences is described by a direct sum of maps as
\begin{equation}\label{pullback}
g^*: H^{q-r-2k}(\tilde Y^{(2k+r+1)},\Q)(-r-k) \to H^{q-r-2k}(\tilde T^{(2k+r+1)},\Q)(-r-k).
\end{equation}

Both spectral sequences degenerate at $E_2 = E_{\infty}$. Keeping track of the multiplicities and  the signs for these pullbacks can be rather hard. Let suppose that locally the defining equations for $Y$ and $T$ are $t = \prod_i z_i^{e_i}$ and $t = \prod_j w_j^{e'_j}$ respectively, and we are given strata $Y_I = Y_{i_1} \cap \ldots \cap Y_{i_p}$ ($i_1 < \ldots < i_p$) and $Y_J = Y_{j_1} \cap \ldots \cap Y_{j_p}$. Then the computation of the multiplicities involved in $g^*: H^*(Y_I,\Q) \to H^*(T_J,\Q)$ essentially amounts to determine the coefficients of $\frac{dw_{j_1}}{w_{j_1}} \wedge \ldots \wedge \frac{dw_{j_p}}{w_{j_p}}$ in $g^*(\frac{dz_{j_1}}{z_{j_1}} \wedge \ldots \wedge \frac{dz_{j_p}}{z_{j_p}})$. This technique will be frequently used in the paper.\vspace{.1in} 

As an example, we describe the map \eqref{pullback} when $f: X \to S$ is a  degeneration of curves with normal crossings singularities on its special fibre $Y$ and $Z$ is the blow-up of $X$ at a closed point $P \in Y$. Let $g: Z \to X$ be the blowing up map. If $P$ is a regular point in the special fibre, the number of components of the special fibre $T$ of $Z$ will simply increase by one (the exceptional divisor $E$) and the remaining components are the same as for $Y$. Hence $g^*: H^0(\tilde Y^{(1)},\C) \to H^0(\tilde T^{(1)},\C)$ is simply the map $g^*(1_{Y_i}) = 1_{T_i} + 1_E$ on the components. 

Let suppose instead that $P$ is singular. Since the description of $g^*$ is local around each closed point, we may assume that the degeneration $f$ is given, in a neighborhood of $P$, by the equation $z_1^{e_1}z_2^{e_2} = t$, being $t$ a chosen parameter on the disk $S$ and $e_1,e_2$ positive integers. Let assume that $e_1 \le e_2$. Then, locally around $P$: $\tilde Y^{(1)} = Y_1 \coprod Y_2$. Set-theoretically one has $Y_i = \{z_i = 0\}$ ($i = 1,2$) and $\tilde Y^{(2)} = Y_1 \cap Y_2 = \{P\}$. Then, $\tilde T^{(1)} = T_1 \coprod T_2 \coprod T_3$ where $T_1$ and $T_2$ are the strict transforms of the two components $Y_i$, while $T_3$ represents the exceptional divisor. We implicitly have fixed the standard orientation on $\tilde Y^{(r)}$ (\eg $\tilde Y^{(2)} = Y_1 \cap Y_2 = Y_{12}$). On $\tilde T^{(r)}$, we choose the orientation for which the exceptional component $T_3$ is always considered as the last one.   

There are only three graded complexes $gr^L_*A^\bu_{\C}$ non zero both on $X$ and $Z$. On $X$ they have the following description
\[
gr^L_{-1}A^\bu_{X,\C} \simeq (a_2)_*\Omega^\bu_{\tilde Y^{(2)}}[-1] 
\]
\[
gr^L_0A^\bu_{X,\C} \simeq (a_1)_*\Omega^\bu_{\tilde Y^{(1)}}
\]

\noindent and via the isomorphism \eqref{nu} one has:
\begin{equation*}
\tilde\nu: gr^L_{1}A^\bu_{X,\C} \overset{\simeq}\to gr^L_{-1}A^\bu_{X,\C}.
\end{equation*}

Hence $E_1^{1,q-1} = \H^q(Y,gr^L_{-1}A^\bu_{X,\C}) = 0$ unless $q = 1$, in which case we get
\[
g^*: H^0(\tilde Y^{(2)},\C) \to H^0(\tilde T^{(2)},\C).
\]

To understand the description of this map, one has to look at the local geometry of the blow-up at $P$. It is quite easy to check that $Z$ is covered by two open sets, say $Z = U \cup V$. To make the notations easier, let call $t_1 = \frac{z_1}{z_2}$ and $t_2 = \frac{z_2}{z_1}$. On $U$, described by $t_2^{e_2} = \frac{t}{z_1^{e_1+e_2}}$, one has coordinates $\{t_2, z_1\}$, $T^{\text{red}}_2 = \{t_2 = 0\}$ and $T^{\text{red}}_3 = \{z_1 = 0\}$. On $V$, described by $t_1^{e_1}  = \frac{t}{z_2^{e_1+e_2}}$, one has coordinates $\{t_1, z_2\}$, $T^{\text{red}}_1 = \{t_1 = 0\}$  and $T^{\text{red}}_3 = \{z_2 = 0\}$. Then $\tilde T^{(2)} = T_{13} \coprod T_{23}$, here we denote $T_{ij} = T_i \cap T_j$. 

On $U$ we have $g^*(\frac{dz_1}{z_1} \wedge \frac{dz_2}{dz_2}) = \frac{dz_1}{dz_1} \wedge \frac{dt_2}{t_2}$, whereas on $V$ one gets $g^*(\frac{dz_1}{z_1} \wedge \frac{dz_2}{dz_2}) = \frac{dt_1}{dt_1} \wedge \frac{dz_2}{z_2}$. Hence, keeping in account the fixed orientation among the components of $T$, the description of the pullback $g^*(1_{\tilde Y^{(2)}}) = g^*(1_{Y_{12}})$ is given by
\[
g^*(1_{Y_{12}}) = 1_{T_{13}} - 1_{T_{23}}.
\]

The presence of a negative sign is due to the change of orientation. This description defines the above map $g^*$ on $H^0$. Similarly, we find that
\[
g^*: H^0(\tilde Y^{(1)},\C) \to H^0(\tilde T^{(1)},\C)
\]

\noindent is given by $g^*(1_{Y_1}) = 1_{T_1} + 1_{T_3}$ and $g^*(1_{Y_2}) =    1_{T_2} + 1_{T_3}$. The description of $g^*$ on the terms $H^1$ goes in parallel.\vspace{.1in}

Let now consider the proper map that $g$ induces on the closed fibres. For simplicity of notations we call it $g: T \to Y$. Let $d = (\dim~T - \dim~Y)$. The above arguments have shown that $g$ induces a pullback map $g^*$ between the cohomologies of the strata: \cf~\eqref{pullback}. Since each stratum is a smooth projective complex variety (not connected), we can use the Poicar\'e duality to associate to each pullback in \eqref{pullback} that contributes to the definition of the map $g^*$ its dual so that we naturally obtain a dual pushforward on the $E_1$--terms of the spectral sequence of weights that is described by a direct sum of maps as
\begin{equation}\label{pushforward}
g_!: H^{q-r-2(k-d)}(\tilde T^{(2k+r+1)},\Q)(-r-k+d) \to H^{q-r-2k}(\tilde Y^{(2k+r+1)},\Q)(-r-k).
\end{equation}\vspace{.1in}

On each stratum $g_!$ is defined by the following formula
\[
(\frac{1}{2\pi\sqrt{-1}})^{d-2k-r}\int_{\tilde Y^{(2k+r+1)}}g_!(\alpha) \cup \beta = (\frac{1}{2\pi\sqrt{-1}})^{2d-2k-r}\int_{\tilde T^{(2k+r+1)}}\alpha  \cup g^*(\beta)
\]\vspace{.1in}

\noindent where $\int$ denotes the morphism trace described by the cap--product with the fundamental class of each component of the stratum, for any chosen couple of elements $\alpha \in H^{q+2(2d-2k-r)}(\tilde T^{(2k+r+1)},\Q(2d-2k-r))$ and $\beta \in H^{-q}(\tilde Y^{(2k+r+1)},\Q)$, $q \in Z,~q \ge 0$.\vspace{.1in}

Notice that although we have a notion of bifiltered pullback 
\[
g^*: (A^\bu_X, L, F) \to (A^\bu_Z, L, F)
\]

\noindent this does not imply a canonical definition of a product structure  on $A^\bu_{\C}$ obtained via pullback along the diagonal map $\Delta: X \to X \times_S X$. In fact, the property of $f: X \to S$ to have normal crossings reduction is not preserved by the product map $f \times f: X \times_S X \to S$. The space $X \times_S X$ is in general not even smooth over $\C$! Finally, we remark that although the monodromy filtration is not multiplicative on the level of the filtered complexes $(A^\bu_{\C},
L)$ (the simple example showed below will motivate this claim), it becomes multiplicative on the limiting cohomology with its mixed Hodge structure. \vspace{.1in}

\begin{ex}\label{ex1}\end{ex}

Let $f: \Pr^1_S \to S$ be a $\Pr^1$-fibration over a disk $S$. We blow a closed point $P \in \Pr^1_0 = Y$ in the fibre $\Pr^1_0$ over the origin $\{0\}$. The resulting map $h: Z \to S$ has a normal crossings special fibre  $h^{-1}(0) = T = T_1 \cup T_2$, where $T_1$ is the strict transform of $Y$ and $T_2$ is the exceptional component (\ie $\Pr^1$). The intersection $Q = T_1 \cap T_2 = T_{12}$ is transverse. Locally around $Q$, $h$ has the following description
\[
h(z_1,z_2) = z_1z_2.
\]

Consider the subcomplex $W_0(A^\bu_{Z,\C})$ of $A^\bu_{Z,\C}$ filtered by the monodromy filtration $L$ induced on it by the one on $A^\bu_{Z,\C}$ (\cf~$\S$~\ref{0},~\eqref{weight}). Its hypercohomology computes $H^*(Y,\C)$ and it can be determined in terms of the homology of the complex 
\begin{equation*}
\{\mathcal C^\bu :~ H^\rnd(\tilde T^{(1)},\C) \overset{d}\to H^\rnd(\tilde T^{(2)},\C)\} \quad =\quad \{\mathcal C^\bu :~H^\rnd(T_1,\C) \oplus H^\rnd(T_2,\C) \overset{d}\to H^\rnd(T_{12},\C)\}
\end{equation*}

\noindent where $\mathcal C^\bu$ sits in degrees zero and one. The differential $d$ on $\mathcal C^\bu$ is of ``\v Cech type'' \ie it is an alternate sum of pullback maps as defined in \eqref{diff}. A product in the filtered derived category $(A^\bu_{Z,\C},L)$ if any exists, should induce a product on $\mathcal C^\bu$. The tensor product $\mathcal C^\bu \otimes  \mathcal C^\bu$ is a complex sitting in degrees zero, one and two and it has the following description
\begin{multline*}
\{\mathcal C^\bu \otimes \mathcal C^\bu :~ \bigoplus_{i,j \in [1,2]}(H^\rnd(T_i,\C)\otimes H^\rnd(T_j,\C)) \overset{d\otimes d}\to \bigoplus_{i=1}^2 \{(H^\rnd(T_i,\C) \otimes H^\rnd(T_{12},\C)) \oplus\\ \oplus (H^\rnd(T_{12},\C) \otimes H^\rnd(T_i,\C))\} \overset{d\otimes d}\to (H^\rnd(T_{12},\C) \otimes H^\rnd(T_{12},\C))\}.
\end{multline*}

There is no canonical way to define a map
\[
\mu: \mathcal C^\bu \otimes \mathcal C^\bu \to \mathcal C^\bu.
\]

In fact, if any existed it should satisfy in degree zero
\begin{align*}
H^\rnd(T_1,\C) \otimes H^\rnd(T_1,\C) &\mapsto H^\rnd(T_1,\C),\qquad H^\rnd(T_2,\C) \otimes H^\rnd(T_2,\C) \mapsto H^\rnd(T_2,\C),\qquad\\\\ &H^\rnd(T_i,\C) \otimes H^\rnd(T_j,\C) \mapsto 0,\quad i,j = 1,2.
\end{align*}

Whereas in degree one would get the non canonical setting
\begin{align*}
H^\rnd(T_1,\C) \otimes &H^\rnd(T_{12},\C) \mapsto H^\rnd(T_{12},\C),\qquad H^\rnd(T_{12},\C) \otimes H^\rnd(T_2,\C) \mapsto H^\rnd(T_{12},\C),\qquad\\\\ 
&H^\rnd(T_2,\C) \otimes H^\rnd(T_{12},\C) \mapsto 0,\qquad H^\rnd(T_{12},\C) \otimes H^\rnd(T_1,\C) \mapsto 0.
\end{align*}\vspace{.1in}

\section{Semistable degenerations with double points}\label{c}

This section is mainly devoted to the determination of $[N]$ for one--dimensional semistable fibrations with at worst double points as singularities. The description of $[N]$ is obtained via the introduction of the algebraic correspondence-type square on the cohomology groups of the special fibre as described in \eqref{2}. A one--dimensional double point degeneration is the simplest example of a normal crossings fibration. The generalization of these results to double points semistable degenerations of arbitrary dimension is done at the end of this paragraph where we also report as an example of application of these results the case of a Lefschetz pencil.\vspace{.1in}

We keep the same notations as in $\S$~\ref{b}, in particular we denote by  $f: X \to S$ a {\it semistable} fibration of fibre dimension one. Its special fibre is denoted by $Y$. Then, locally around a double point $P \in Y$ the description of $f$ looks like
\[
f(z_1,z_2) = z_1z_2
\]

\noindent for $\{z_1,z_2\}$ local parameters on $X$ at $P$. 
For one dimensional fiberings, the only group where the local monodromy may act non trivially is $gr^L_2H^1(\tsX,\Q)$, in which case the identity map on the $E_1$-terms of the weight spectral sequence \eqref{1a}
\begin{equation*}
E_1^{-1,2} = H^0(\tilde Y^{(2)},\Q)(-1) \overset{\text{Id}}\to H^0(\tilde Y^{(2)},\Q)(-1) = E_1^{1,0}(-1)
\end{equation*}

\noindent determines an isomorphism of rational Hodge structures of weight two on the related graded groups $E_2 = E_{\infty}$. This isomorphism is induced by the action of the local monodromy $N$ around the origin: 
\begin{equation*}
N: gr^L_2H^1(\tsX,\Q) \overset{\simeq}\to (gr^L_0H^1(\tsX,\Q))(-1) 
\end{equation*}

It is a well known consequence of the Clemens-Schmid exact sequence (considered as a sequence of mixed Hodge structures) that 
\[
gr^L_2H^1(\tsX,\Q) \neq 0~\Leftrightarrow~\text{Ker}(\rho^{(2)}: H^1(\tilde Y^{(1)},\Q) \to H^1(\tilde Y^{(2)},\Q)) \neq 0~\Leftrightarrow~h^1(|\Gamma|) \neq 0
\]

\noindent where $h^1(|\Gamma|)$ is the dimension of the first cohomology group of the geometric realization of the dual graph of $Y$ (\cf~\cite{M}).
It follows from proposition~\ref{prop1} that $[N] \in \H^2(T,(gr^L_{-2}A^\bu_{Z,\Q})) = H^{0}(\tilde{T}^{(3)},\Q)$ determines a Hodge class  
\begin{multline}\label{N}
[N] \in Hom_{MH}(\Q(0),gr^L_0H^2(T,\Q)) \simeq \\\simeq Hom_{MH}\biggl (\Q(0),\frac{H^{0}(\tilde{T}^{(3)},\Q)}{\text{Image}(\rho^{(3)}: H^{0}(\tilde{T}^{(2)},\Q) \to H^{0}(\tilde{T}^{(3)},\Q))}\biggr ).
\end{multline}

Here $T$ is the special fibre of a normal--crossings degeneration $h: Z \to S$. The variety $Z$ is a smooth threefold over $\C$ obtained via resolution of the singularities of $X \times_S X$. Notice that no more 
than three components of $T$ intersect at the same closed point since $\dim Z = 3$.

We shall determine the Hodge cycle $[N] \in E_1^{2,0}(Z) = H^{0}(\tilde{T}^{(3)},\Q)$ by means of a ``correspondence type'' map 
\[
N: \H^*(Y, gr^L_rA^\bu_{X,\Q}) \to \H^*(Y, (gr^L_{r-2}A^\bu_{X,\Q})(-1)) = \H^*(Y, gr^L_r(A^\bu_{X,\Q}(-1)))
\]

\noindent as we explained in \eqref{2}. We will prove that the map $N$ is zero for $* \neq 1$ and is the identity for $* = 1 = r$. On the $E_2$-level it will induce (for $* = 1 = r$) a commutative diagram 
\begin{equation*}
\begin{CD}
gr^L_2H^1(\tsX \times \tsX, \Q) @>{[N]~\cdot}>> gr^L_2H^{3}(\tsX \times \tsX, \Q) = E_2^{1,2} \\
@A{(p_1)^*}AA @VV{(p_2)_*}V \\
E_2^{-1,2} = gr^L_2H^1(\tsX, \Q) @>{N}>> (gr^L_0H^1(\tsX, \Q))(-1) = (E_2^{1,0})(-1)
\end{CD}
\end{equation*}

The pullback $p_1^*$ and pushforward $(p_2)_*$ are defined as in $\S$~\ref{b}. The above diagram will determine uniquely both $[N] \in Hom_{MH}(\Q(0),gr^L_0H^2(T,\Q))$ and the product $[N]~\cdot$.\vspace{.1in}

The following result defines the geometry of the model $Z$ and the special fibre $T$ after resolving the singularities of $X \times_S X$ and $Y \times Y$.

\begin{lem}\label{lem0} Let $z_1z_2 = w_1w_2$ be a local description of $X \times_S X$ around the point $(P,P)$, with $P \in Y = Y_1 \cup Y_2$ a double point of $f$ and $\{w_1,~w_2\}$ a second set of regular parameters on $X$ at $P$. After a blow-up of $X \times_S X$ with center at the origin $(z_1,z_2,w_1,w_2)$, the resulting degeneration $h: Z \to S$ is normal--crossings. Its special fibre $T$ is the union of five irreducible components: $T = \cup_{i=1}^5 T_i$. We number them so that the first four are the strict transforms of the irreducible components $Y_i \times Y_j$ of $Y \times Y$, namely $T_1 = (Y_1 \times Y_1)^{\tilde{}}$, $T_2 = (Y_1 \times Y_2)^{\tilde{}}$, $T_3 = (Y_2 \times Y_1)^{\tilde{}}$, $T_4 = (Y_2 \times Y_2)^{\tilde{}}$. The last one $T_5$ represents the exceptional divisor of the blow-up. We have $\tilde T^{(1)} = \coprod_i T_i$. The scheme $Z$ is covered by four affine charts $\mathcal U_j$. On each of them there are three  non empty components $T_k$. The scheme $\tilde T^{(3)}$ is made by the disjoint union of four zero dimensional schemes (closed points): $T_{125} \in \mathcal U_2$, $T_{135} \in \mathcal U_4$, $T_{245} \in \mathcal U_3$ and $T_{345} \in \mathcal U_1$, each of whose supports projects isomorphically onto the diagonal $\Delta_{12}: Y_{12} \to Y_{12} \times Y_{12}$. 
\end{lem}
\begin{proof} The local description of $X \times_S X$ around $(P,P)$ is given by the equations $z_1z_2 = w_1w_2$ and $z_1z_2 = t$, for $t \in S$ a fixed parameter on the disk. We choose the standard orientation of the sets $\{z_1,z_2\}$ and $\{w_1,w_2\}$ and we write $w_{i_1}' = \frac{w_i}{z_1}$, $w_{i_2}' = \frac{w_i}{z_2}$, $w_{ij} = \frac{w_i}{w_j}$, $z_{i_1}' = \frac{z_i}{w_1}$, $z_{i_2}' = \frac{z_i}{w_2}$ and $z_{ij} = \frac{z_i}{z_j}$, for $i,j = 1,2$. After a single blow-up of $X \times_S X$ at the origin $(z_1,z_2,w_1,w_2)$, the resulting model $Z$ is non singular as one can see by looking at the first of the following tables which describes $Z$ on each of the four charts $\mathcal U_j$ who cover it. In the second table, we have collected for each $\mathcal U_j$, the description of the non empty divisors $T_k \in T^{(1)}$ there. We use the pullbacks $p_1^*(\frac{dz_1}{z_1} \wedge \frac{dz_2}{z_2})$ and $p_2^*(\frac{dw_1}{w_1} \wedge \frac{dw_2}{w_2})$ to define in the third table the pullbacks $p_i^*(1_{Y_{12}}) \in H^0(\tilde T^{(2)},\Q)$.
\begin{center}
\begin{tabular}{c|c} 
Open sets & Loc. coordinates and relations \\ \hline\hline
  {\rule[-3mm]{0mm}{8mm}} 
$\mathcal U_1$ & $\{w_{1_1}',w_{2_1}',z_1\},~w_{1_1}'w_{2_1}' = z_{21}$\\ \hline
  {\rule[-3mm]{0mm}{8mm}} 
$\mathcal U_2$ &$\{w_{1_2}',w_{2_2}',z_2\},~w_{1_2}'w_{2_2}' = z_{12}$  \\ \hline
  {\rule[-3mm]{0mm}{8mm}} 
$\mathcal U_3$ &$\{z_{1_1}',z_{2_1}',w_1\},~z_{1_1}'z_{2_1}' = w_{21}$
  \\ \hline
  {\rule[-3mm]{0mm}{8mm}} 
$\mathcal U_4$ &$\{z_{1_2}',z_{2_2}',w_2\},~z_{1_2}'z_{2_2}' = w_{12}$  \\ \hline
\end{tabular}
\end{center}\vspace{.2in}

\begin{center}
\begin{tabular}{c|c} 
Open sets & 
 Divisors\\ \hline\hline
  {\rule[-3mm]{0mm}{8mm}} 
$\mathcal U_1$ & $T_3 = \{w_{1_1}' = 0\},~T_4 = \{w_{2_1}' = 0\},~T_5 = \{z_1 = 0\}$   \\ \hline
  {\rule[-3mm]{0mm}{8mm}} 
$\mathcal U_2$ & $T_1 = \{w_{1_2}' = 0\},~T_2 = \{w_{2_2}' = 0\},~T_5 = \{z_2 = 0\}$ \\ \hline
  {\rule[-3mm]{0mm}{8mm}} 
$\mathcal U_3$ & $T_2 = \{z_{1_1}' = 0\},~T_4 = \{z_{2_1}' = 0\},~T_5 = \{w_1 = 0\}$ \\ \hline
  {\rule[-3mm]{0mm}{8mm}} 
$\mathcal U_4$ & $T_1 = \{z_{1_2}' = 0\},~T_3 = \{z_{2_2}' = 0\},~T_5 = \{w_2 = 0\}$ \\ \hline
\end{tabular}
\end{center}\vspace{.2in}

\begin{center}
\begin{tabular}{c|c|c} 
Open sets & $p_1^*(1_{Y_{12}})$ & $p_2^*(1_{Y_{12}})$\\ \hline\hline
  {\rule[-3mm]{0mm}{8mm}} 
$\mathcal U_1$ & $-1_{T_{35}}-1_{T_{45}}$ & $-1_{T_{45}}+1_{T_{35}}+1_{T_{34}}$ \\ \hline
  {\rule[-3mm]{0mm}{8mm}} 
$\mathcal U_2$ & $1_{T_{15}}+1_{T_{25}}$ & $-1_{T_{25}}+ 1_{T_{15}}+1_{T_{12}}$  \\ \hline
  {\rule[-3mm]{0mm}{8mm}} 
$\mathcal U_3$ & $-1_{T_{45}}+1_{T_{25}}+1_{T_{24}}$ & $-1_{T_{25}}-1_{T_{45}}$ \\ \hline
  {\rule[-3mm]{0mm}{8mm}} 
$\mathcal U_4$ & $-1_{T_{35}}+1_{T_{15}}+1_{T_{13}}$ & $1_{T_{15}}+1_{T_{35}}$ \\ \hline
 \end{tabular}
\end{center}\vspace{.1in}

The global description of the pullbacks $p_1^*(1_{Y_{12}})$ and $p_2^*(1_{Y_{12}})$ is
\[
p_1^*(1_{Y_{12}}) = (1_{T_{15}} + 1_{T_{25}} - 1_{T_{35}} - 1_{T_{45}}) + 1_{T_{13}} + 1_{T_{24}}
\]
\[
p_2^*(1_{Y_{12}}) = (1_{T_{15}} - 1_{T_{25}} + 1_{T_{35}} - 1_{T_{45}}) + 1_{T_{12}} + 1_{T_{34}}.
\]

Finally, notice that each $\mathcal U_j$ is isomorphic to $\Af^3$ and in each of them one has three non empty components $T_k$.
\end{proof}\vspace{.1in}

The following result holds
\begin{thm}\label{th0} Let $f: X \to S$ be the semistable degeneration of curves as described above. Then, the following description of $[N] \in H^0(\tilde T^{(3)},\Q)$ (\cf~\eqref{N}) holds:
\[
[N] = a_{125}1_{T_{125}} + a_{135}1_{T_{135}} + a_{245}1_{T_{245}} + a_{345}1_{T_{345}}
\]

\noindent where the (rational) numbers a's are subject to the following requirement:
\[
-2a_{125} + 2a_{135} - 2a_{245} + 2a_{345} = 1.
\]

The induced class $[N]$ in $gr^L_0H^2(T,\Q)$ (\ie modulo boundary relations via the restriction map $\rho^{(3)}$ \cf~\eqref{diff}) determines a unique zero--cycle.
\end{thm}

\begin{proof} We determine $[N]$ as a cocycle making the following square commute 
\begin{equation}\label{dia2}
\begin{CD}
gr^L_2H^1(\tsX \times \tsX, \Q) @>{[N]~\cdot}>> gr^L_2H^{3}(\tsX \times \tsX, \Q) = E_2^{1,2} \\
@A{(p_1)^*}AA @VV{(p_2)_*}V \\
E_2^{-1,2} = gr^L_2H^1(\tsX, \Q) @>{N}>> (gr^L_0H^1(\tsX, \Q))(-1) = (E_2^{1,0})(-1)
\end{CD}
\end{equation}

In terms of cohomologies of strata, we have to describe explicitly a representative of $[N]$ in  $E_1^{2,0}(Z)$ that satisfies the commutativity of 
\begin{equation}\label{dia1}
\begin{CD}
H^0(\tilde T^{(2)},\Q)(-1) @>{[N]~\cdot}>> H^{2}(\tilde T^{(2)},\Q) \\
@A{p_1^*}AA @VV{(p_2)_*}V \\
E_1^{-1,2} = H^0(\tilde Y^{(2)},\Q)(-1) @= H^0(\tilde Y^{(2)},\Q)(-1) = E_1^{1,0}(-1).
\end{CD}
\end{equation}

With the notations used in lemma~\ref{lem0} the description of $[N]$ is given by
\[
[N] = a_{125}1_{T_{125}} + a_{135}1_{T_{135}} + a_{245}1_{T_{245}} +  a_{345}1_{T_{345}}.
\]

For the standard choice of the orientations of $\{z_1,z_2\}$ and $\{w_1,w_2\}$ and the numbering of the $T_i$'s defined in lemma~\ref{lem0}, the local description of the pullbacks $p_i^*(1_{Y_{12}})$ for $i = 1,2$ is given in the third table of the above lemma. Following the definition  described in the Appendix (\cf~\eqref{product}), the product $[N] \cdot p_1^*(1_{Y_{12}}(-1))$ is then the following
\begin{equation}\label{pr}
[N] \cdot p_1^*(1_{Y_{12}}(-1)) =
\end{equation}

$$
= [N] \cdot (1_{T_{15}}(-1) + 1_{T_{15}}(-1) - 1_{T_{35}}(-1) - 1_{T_{45}}(-1)) =
$$
\begin{align*}
\begin{split}
&= a_{125}(g_1(1_{T_{125}} \cdot 1_{T_{15}}(-1)) - g_2(1_{T_{125}} \cdot 1_{T_{25}}(-1))) + a_{135}(g_1(1_{T_{135}} \cdot 1_{T_{15}}(-1)) ~+\\
&+ g_3(1_{T_{135}} \cdot 1_{T_{35}}(-1))) +  
 a_{245}(g_2(1_{T_{245}} \cdot 1_{T_{25}}(-1)) + g_4(1_{T_{245}} \cdot 1_{T_{45}}(-1)))~+\\
&+ a_{345}(- g_3(1_{T_{345}} \cdot 1_{T_{35}}(-1)) + g_4(1_{T_{345}} \cdot 1_{T_{45}}(-1))) =\\ 
&= a_{125}(1_{T_{25}} - 1_{T_{15}}) + a_{135}(1_{T_{35}} + 1_{T_{15}}) + a_{245}(1_{T_{45}} + 1_{T_{25}}) + a_{345}(- 1_{T_{45}} + 1_{T_{35}}).  
\end{split}
\end{align*}

The maps $g_1$, $g_2$, $g_3$ and $g_4$ are the pushforwards as introduced in  the Appendix. The following formula illustrates the  product $1_{T_{ijk}} \cdot \sum_{l,m}1_{T_{lm}}(-1)$ following the definition of it given in the Appendix: 
\begin{multline*}
1_{T_{ijk}} \cdot \sum_{l,m}1_{T_{lm}}(-1) = 1_{T_{ijk}} \cdot (1_{T_{ik}}(-1) + 1_{T_{jk}}(-1)) = g_i(1_{T_{ijk}} \cdot 1_{T_{ik}}(-1)) - g_j(1_{T_{ijk}} \cdot\\ \cdot 1_{T_{jk}}(-1)) = g_i(1_{T_{ijk}}) - g_j(1_{T_{ijk}}) \in \text{Image}\biggl (\bigoplus_t g_t: H^0(\tilde T^{(3)},\Q)(-1) \to H^2(\tilde T^{(2)},\Q)\biggr ).
\end{multline*}

In \eqref{pr}, we have denoted, for simplicity of notations, the difference $g_i(1_{T_{ijk}}) - g_j(1_{T_{ijk}})$ with $1_{T_{jk}} - 1_{T_{ik}}$. The map $g_i$ represents the pushforward on cycles deduced from the embedding $g_i: T_{ijk} \to T_{jk}$. The definition of $g_j$ is similar.

Therefore, via the local definition of the pushforward $(p_2)_*$ along the affine charts (\cf~$\S\ref{b}$ and third table in lemma~\ref{lem0}), we obtain:
\[
(p_2)_*([N] \cdot p_1^*(1_{Y_{12}}(-1))) = (- 2a_{125} + 2a_{135} - 2a_{245} + 2a_{345})1_{Y_{12}}(-1).
\]

The commutativity of \eqref{dia1} and hence of \eqref{dia2} is then equivalent to the requirement
\[
- 2a_{125} + 2a_{135} - 2a_{245} + 2a_{345} = 1.
\]

Hence, the operator $[N]$ is determined as a cocycle in  $H^0(\tilde T^{(3)},\Q)$ by the setting
\begin{equation}\label{*}
\begin{split}
[N] =~&a_{125}1_{T_{125}} + a_{135}1_{T_{135}} + a_{245}1_{T_{245}} + a_{345}1_{T_{345}};\\\\
 &- 2a_{125} + 2a_{135} - 2a_{245} + 2a_{345} = 1.
\end{split}
\end{equation}

Up to boundary relations by means of the restriction map $\rho^{(3)}$ which connects the elements $1_{T_{125}}$ with $1_{T_{245}}$ and $1_{T_{135}}$ with $1_{T_{345}}$, \eqref{*} determines a unique zero--cycle in the quotient $E_2^{2,0}(Z)$ (\cf~\eqref{N}). Of course, if $N = 0$, this class may be trivial.\end{proof}\vspace{.2in}

\begin{rem}\label{rem1}\end{rem}

The description of $[N] \in E_1^{2,0}(Z)$ as well as the relation among the coefficients $a_{ijk}$ in \eqref{*} is not unique in $E_1$. In fact, it depends on the choice of the desingularization process, as well as on the ordering of the components $T_k \in \tilde{T}^{(1)}$. For example, for the ordering of $T_k$ for which $T_1$ represents in each chart the exceptional divisor of the blow-up ($T_2 = (Y_1 \times Y_1)^{\tilde{}}$, $T_3 = (Y_1 \times Y_2)^{\tilde{}}$, $T_4 = (Y_2 \times Y_1)^{\tilde{}}$, $T_5 = (Y_2 \times Y_2)^{\tilde{}}$), the setting \eqref{*} becomes
\begin{equation*}
\begin{split}
[N] =~&a_{123}1_{T_{123}} + a_{124}1_{T_{124}} + a_{135}1_{T_{135}} + a_{145}1_{T_{145}};\\\\
 &- a_{123} + a_{124} - a_{135} + a_{145} = 1.
\end{split}
\end{equation*}\vspace{.1in}

If instead we choose to desingularize $X \times_S X$ via a blowing-up along $z_1 = w_1 = 0$ and we set the order among the $T_k$'s so that the  exceptional divisor is represented in each chart by the last component (\ie $T_1 = (Y_1 \times Y_2)^{\tilde{}}$, $T_2 = (Y_2 \times Y_1)^{\tilde{}}$, $T_3 = (Y_2 \times Y_2)^{\tilde{}}$, $T_4 = (Y_1 \times Y_1)^{\tilde{}}$,), then we would get 
\begin{equation*}
\begin{split}
[N] =~&a_{134}1_{T_{134}} + a_{234}1_{T_{234}};\\\\
 &- a_{134} + a_{234} = 1.
\end{split}
\end{equation*}

It is a consequence of  the uniqueness of the product structure on the corresponding $E_2$--terms that all these different settings determine a unique description of $[N] \in E_2^{2,0}(Z)$.\vspace{.2in}

The description of $[N]$ for a double points degeneration of higher fibre dimension can be easily deduced from the case worked out for curves. In fact, it was already remarked that the description of $[N]$ in the cohomology of the strata of the special fibre of the fibre product resolution is of local nature, \ie it can be described locally around each double point. For a higher dimensional double points degeneration $[N]$ is again described in terms of a ``diagonal'' cocycle whose support projects isomorphically onto the diagonal $\Delta_{12} \in Y_{12} \times Y_{12}$ as  was shown in theorem~\ref{th0}. In general, that ``diagonal'' cocycle is formally locally a bundle over the corresponding diagonal cocycle which comes up for a degeneration of curves. This is a consequence of the local description of the degeneration map around a double point. We give some details for these claims.\vspace{.1in}

Let $f: X \to S$ be a semistable degeneration with double points of fibre dimension $d$ over the disk $S$. Then, locally in a neighborhood of a double point $P$ on $Y$, $f$ has the following description
\[
f(z_1,\ldots,z_n) = z_iz_j
\]

\noindent for $\{z_1,\ldots,z_n\}$ a set of regular parameters on $X$ at $P$ and suitable indices $i < j$ in $I = \{1,\ldots,n\}$. Let $Y = Y_1 \cup Y_2$ be the local description of $Y$ in a neighborhood of $P \in Y_1 \cap Y_2 = Y_{12}$. Because $\{z_1,\ldots,\hat{z}_i,\ldots,\hat z_j,\ldots, z_n\}$ are free parameters for this description, formally locally around $P$, the special fibre is isomorphic to $\Af^{d-1} \times \hat Y$ with $\hat Y = \hat{Y_1} \cup \hat{Y_2}$ of dimension $1$. In terms of local coordinates, $Y$ is described as $\text{Spec}\biggl (\C\{\{z_1,\dots,\hat{z_i},\ldots,\hat{z_j},\ldots z_n\}\}[z_i,z_j]/z_iz_j\biggr )$. The model $X$ is formally locally isomorphic to $\Af^{d-1} \times \hat X$, with $\hat X$ of fibre dimension $1$ and special fibre $\hat Y$. The formal description of $X \times_S X$ is similar, namely $X \times_S X \simeq \Af^{d-1} \times \Af^{d-1} \times (\hat X \times_S \hat X)$. Keeping  the same notations introduced before, we get a formal local description of the stratum $\tilde T^{(3)}$ containing the cocycle $[N]$ as $\Delta_{\Af^{d-1}} \times \hat{\tilde T}^{(3)}$, with  $\hat{\tilde T}^{(3)}$ made by a collection of points. This scheme maps again isomorphically onto $\Delta_{\hat{Y}_{12}}$, the diagonal in $\hat{Y}_{12} \times \hat{Y}_{12}$.

In this way, the description of $[N]$ can be deduced from a formal local description of the Lefschetz pencil of fibre dimension one $f: \hat X \to \C\{\{t\}\}$.  Hence, we get a formal local class representative of $N$ as a bundle over the diagonal cocycle which describes $[N]$ in theorem~\ref{th0}.  In particular this proves the following

\begin{cor} Let $f: X \to S$ be a semistable double points degeneration of fibre dimension $d$. Then
\begin{equation*}
[N] \in 
\frac{CH^{d-1}(\tilde T^{(3)})}{\text{Image}\biggl (\rho^{(3)}: CH^{d-1}(\tilde{T}^{(2)}) \to CH^{d-1}(\tilde{T}^{(3)})\biggr )}
\end{equation*}
is represented by a unique algebraic cocycle of codimension d-1 in the stratum $\tilde T^{(3)}$.\label{cor1}
\end{cor}\vspace{.1in}

Notice that for a double point degeneration of fibre dimension $d > 1$, $[N]$ may represent the monodromy map acting non trivially on  different  graded pieces of the limiting cohomology. However, they are all of type $gr^L_{q+1}H^q(\tsX,\Q) = E_2^{-1,q+1}(X)$ for $q \in [0,d]$. In fact, for double point degenerations we have always $N = 0$ on $gr^L_qH^q(\tsX,\Q)$, and $\H^*(Y,gr^L_iA^\bu_{X,\C}) = 0$ for $i \neq -1,~0,~1$ because no more than two components of $Y$ intersect simultaneusly at the same closed point.\vspace{.2in}

As an example of application of these results we consider the case of a  Lefschetz pencil of fibre dimension at least three. The description of $[N]$ is the same to the one just described for a degeneration with double points. We will only show how to reduce in this case the study of $[N]$ to the previous one. A Lefschetz pencil of fibre dimension greater than one is not even normal-crossings because the special fibre is irreducible and singular. We will only consider the case of odd fibre dimension since Lefschetz pencils of even fibre dimension have trivial monodromy always. 

Let $f': \mathcal X \to S$ be such a pencil and let $n = 2m + 1$ be the dimension of its fibre. Locally, in a neighborhood of the singular point of the special fibre $\mathcal Y$, the pencil $f'$ is described by
\[
f'(z_0,\ldots,z_n) = \sum_{\nu = 0}^m z_\nu z_{\nu +1+m}
\]

\noindent where as usual $\{z_0,\ldots,z_n\}$ represents a set of regular parameters on $\mathcal X$. It is clear from the definition that the special fibre $\mathcal Y$ is irreducible and singular at the origin $(z_0,\ldots,z_n)$. However, after a single blow-up at that point we get a normal-crossings degeneration $f: X \to S$ with special fibre locally described by $Y = Y_1 \cup Y_2$. The component $Y_1$ is the exceptional divisor of the blow-up, a projective space of dimension $n$ which intersects the strict transform $Y_2$ of $Y$ along a quadric hypersurface $Y_{12}$ of dimension $2m$. The component $Y_1$ appears with multiplicity $e_1 = 2$ whereas $Y_2$ is reduced (\ie $e_2 = 1$). 
Let $h: X \to \mathcal X$ be the blow-up map. It is a (proper) map of $S$-schemes, therefore it induces a morphism
\[
g^*\R\P_{f'}(\Q_{\mathcal X}) \to \R\P_{f}(\Q_{X})
\]

\noindent of complexes of nearby cycles. This morphism induces in turn a homomorphism between the corresponding hypercohomologies
\[
g^*: \H^i(\mathcal Y,\R\P_{f'}(\Q)) \to \H^i(Y,\R\P_{f}(\Q))
\]

In order to work with the resolution $A^\bu_{\Q}$ of $\R\P_{f}(\Q)$ which carries the monodromy filtration, we have to consider $Y$ with its reduced structure (the exceptional divisor has multiplicity $e_1 = 2$ as algebraic cycle on $X$). Because $g.c.d.(e_1,e_2) = 1$ $\forall y \in Y$  the action of the local monodromy on the complex of sheaves $\R\P_{f}(Y,\Q)$ is unipotent (\cf~$\S$~\ref{0}). That implies that the monodromy operator acts unipotently on cohomology.

Because $f'$ is a Lefschetz pencil of fibre dimension $n$, the only group where $N$ acts non trivially is $H^{n}(\tsX,\Q)$. Also, $[N]$ determines an element in $(H^{2n}(\tsX \times \tsX,\Q(n-1)))^{\pi_1}$ and because the generic fibres of $f'$ and $f$ are the same, we may as well consider $[N] \in \H^{2n}(Y\times Y,\R\P_{f}(\Q))^{\pi_1}$.\vspace{.1in} 

The map $f$ is locally described by $z_i^2q(z_0,\ldots,\hat z_i,\ldots,z_n) = t$ for some $i \in [0,n]$, $t$ being a local parameter on $S$ and $q(z_0,\ldots,\hat z_i,\ldots,z_n)$ an irreducible quadratic polynomial. Via the extension of the basis $S' \to S$ $\tau \mapsto \sqrt t$, the  degeneration $f$ is deformed to $w_iz_i = \tau$, with $w_i = \frac{\tau}{z_i}$ and $w_i^2 = h$. It is clear that this procedure does not affect the special fibres (\ie the reduced closed fibres are the same). Hence, after a possible normalization of the resulting model, we obtain a double point semistable degeneration $h: Z \to S$. Let $T = T_1 \cup T_2$ be its special fibre. Then $[N]$ can be seen as a Hodge cycle in $H^{2n}(T\times T,\R\P_{h}(\Q))^{\pi_1} = \text{Ker}(\tilde N) \cap H^{2n}(\tsX \times \tsX,\Q)$, for $\tilde N = 1 \otimes N + N \otimes 1$. The geometric description of $[N]$ is then the same as the one we have shown before. The class $[N]$ represents the monodromy operator acting non trivially only on $gr^L_{n+1}H^n(\tsX,\Q)$. \vspace{.2in}

\section{Semistable degenerations with triple points}\label{e}

A semistable degeneration with triple points is the first case where both the operators $N$ and $N^2$ may be non trivial. In this paragraph we will mainly consider a triple point degeneration of surfaces. The description of $[N]$ and $[N^2]$ for higher dimensional triple points degenerations can be deduced from the one for surfaces using the same kind of arguments described in the last paragraph for double points degenerations of higher fibre dimension.\vspace{.1in}

Let $f: X \to S$ be a surfaces degeneration with reduced normal crossings and with triple points on its special fibre $Y$. We keep the basic notations as in the previous sections. Then, locally around a triple point $P \in Y$ we may assume that $f$ has the following description:
\[
f(z_1,z_2,z_3) = z_1z_2z_3.
\]

As usual, $\{z_1,z_2,z_3\}$ is a regular set of parameters on $X$ at $P$. Globally on $X$, the special fibre can be the union of more than three components \ie $Y = Y_1 \cup\ldots\cup Y_N$, but at most three of them  intersect at the same closed point. The Clemens--Schmid exact sequence of mixed Hodge structures  describes the behavior of the operators $N$ and $N^2$ in terms of some invariants on the special fibre. Namely
\begin{lem} (Monodromy criteria) Let $f: X \to S$ be a semistable degeneration of surfaces, then
\begin{align*}
N = 0\qquad\text{on}\quad H^1(\tsX,\Q) & \Leftrightarrow h^1(|\Gamma|) = 0\\
N = 0\qquad\text{on}\quad H^2(\tsX,\Q) & \Leftrightarrow h^2(|\Gamma|) = 0\qquad\text{and}\quad \rho^{(2)}: H^1(\tilde Y^{(1)},\Q) \twoheadrightarrow H^1(\tilde Y^{(2)},\Q)\\
N^2 = 0\qquad\text{on}\quad H^2(\tsX,\Q) &\Leftrightarrow h^2(|\Gamma|) = 0.
\end{align*}

Here $h^i(|\Gamma|)$ denotes the dimension of the ith-cohomology group of the geometric realization of the dual graph of $Y$. 
\end{lem}
\begin{proof} \cf~\cite{M}.\end{proof}

A degeneration of K-3 surfaces with special fibre made by rational surfaces intersecting along a cycle of rational curves, is an example  for which both $N$ and $N^2$ are non zero (\cf~\cite{PP}). 

Let suppose that at least one between $gr^L_2H^1(\tsX,Q)$ and $gr^L_3H^2(\tsX,\Q)$ is non zero (for the above example it is well known that $gr^L_2H^1(\tsX,\Q) = 0$, as $H^1(\tsX,\Q) = 0$). The map $N$ acts on them as an isomorphism of pure Hodge structures 
\begin{equation*}
N: gr^L_2H^1(\tsX,Q) \overset{\simeq}\to (gr^L_0H^1(\tsX,\Q))(-1)
\end{equation*}
\begin{equation*}
N: gr^L_3H^2(\tsX,\Q) \overset{\simeq}\to (gr^L_1H^2(\tsX,\Q))(-1).
\end{equation*}

The only group where $N^2$ behaves as an isomorphism is $gr^L_4H^2(\tsX,\Q)$. The map $N^2$ is defined by the composition
\begin{equation*}
gr^L_4H^2(\tsX,\Q) \overset{N}\to (gr^L_2H^2(\tsX,\Q))(-1) \overset{N}\to (gr^L_0H^2(\tsX,\Q))(-2).
\end{equation*}

 The sequence is not exact in the middle. The map $N$  on the left is injective and the one on the right surjects $(gr^L_2H^2(\tsX,\Q))(-1)$ onto $(gr^L_0H^2(\tsX,\Q))(-2)$. Its kernel, in term of the spectral sequence of weights is  
$$
\biggl (\text {Im}(\H^2(Y,gr^W_1\Omega_X^{\bu +1}(\log Y)) \otimes \Q \to \H^2(Y,A^\bu_{X,\Q}))\biggr )(-1) \simeq 
$$
$$
\simeq \frac{\text{Ker}(\rho^{(2)}: H^2(\tilde Y^{(1)},\Q)(-1) \to H^2(\tilde Y^{(2)},\Q)(-1))}{\text{Im}(\gamma^{(2)}: H^0(\tilde Y^{(2)},\Q)(-2) \to H^2(\tilde Y^{(1)},\Q)(-1))}.
$$\vspace{.1in}

We first consider $N$ and its related class $[N]$. Both $gr^L_2H^1(\tsX,\Q)$ and $gr^L_3H^2(\tsX,\Q)$ are described in terms of cohomology classes on $\tilde Y^{(2)}$ (\cf~\eqref{1b}). The study of the correspondence-diagram \eqref{2} is similar for them. Namely, once one has found an algebraic  cycle representing $[N]$, it certainly makes both the correspondence diagrams commute. For degenerations of surfaces it follows from proposition~\ref{prop1} that
\begin{equation}\label{4}
[N] \in  (gr^L_2H^4(T,\Q))(1) \simeq
\frac{\text{Ker}(\rho^{(4)}: H^2(\tilde T^{(3)},\Q)(1) \to H^2(\tilde T^{(4)},\Q)(1))}{\text{Im}(\rho^{(3)}: H^2(\tilde T^{(2)},\Q)(1) \to H^2(\tilde T^{(3)},\Q)(1))}
\end{equation}

\noindent where $h: Z \to S$ is a normal crossings degeneration with special fibre $T$ and generic fibre $\tsX \times \tsX$ obtained via resolution of the singularities of $X \times_S X$. Similarly, one has
\begin{equation}\label{5}
[N^2] \in  gr^L_0H^4(T,\Q) \simeq
\frac{H^0(\tilde T^{(5)},\Q)}{\text{Im}(\rho^{(5)}: H^0(\tilde T^{(4)},\Q) \to H^0(\tilde T^{(5)},\Q))}.
\end{equation}\vspace{.1in}

Both $[N]$ and $[N^2]$ have the further property to be Hodge cycles in the  cohomologies of the corresponding strata. The following lemma determines the geometry of the model $Z$ and the special fibre $T$ after resolving the singularities of $X \times_S X$ and $Y \times Y$.

\begin{lem}\label{lem1} Let $z_1z_2z_3 = w_1w_2w_3$ be a local description of $X \times_S X$ around the point $(P,P)$, being $P \in Y = \cup_{i=1}^3 Y_i$ a triple point of $f$ and $\{w_1,~w_2,~w_3\}$ a second set of regular parameters on $X$ at $P$. After three blows-up of $X \times_S X$ with centers at $z_i = 0 = w_i$ ($i = 1, 2, 3$) the resulting degeneration $h: Z \to S$ is normal--crossings. Its special fibre $T$ is the union of nine irreducible components: $T = \cup_{i=1}^9 T_i$. We number them so that the first six are the strict transforms of the irreducible components $Y_i \times Y_j$ of $Y \times Y$: $T_1 = (Y_1 \times Y_2)^{\tilde{}}$, $T_2 = (Y_1 \times Y_3)^{\tilde{}}$, $T_3 = (Y_2 \times Y_1)^{\tilde{}}$, $T_4 = (Y_2 \times Y_3)^{\tilde{}}$, $T_5 = (Y_3 \times Y_1)^{\tilde{}}$, $T_6 = (Y_3 \times Y_2)^{\tilde{}}$. The last three components are the exceptional divisors of the three blows-up: $T_7 = (Y_1 \times Y_1)^{\tilde{}}$, $T_8 = (Y_2 \times Y_2)^{\tilde{}}$, $T_9 = (Y_3 \times Y_3)^{\tilde{}}$ . We have $\tilde T^{(1)} = \coprod_i T_i$. The scheme $Z$ is covered by eight affine charts, on each of them there are at most five  non empty components $T_i$. Among the components $T_{ijk}$ whose disjoint union defines the scheme $\tilde T^{(3)}$, $T_{178}$ and $T_{378}$ contain \resp the curves ``diagonal'' $\tilde{\delta}_{12}$ and $\delta_{12}$ whose supports project isomorphically onto the diagonal $\Delta_{12}: Y_{12} \to Y_{12} \times Y_{12}$. Similarly, $T_{279}$ and $T_{579}$ contain \resp $\tilde{\delta}_{13}$ and $\delta_{13}$ whose support projects isomorphically onto $\Delta_{13}: Y_{13} \to Y_{13} \times Y_{13}$. Finally, $T_{489}$ and $T_{689}$ contain $\tilde{\delta}_{23}$ and $\delta_{23}$ whose support is isomorphic to $\Delta_{23}$. The exceptional surface $T_{789}$--intersection of the three exceptional divisors of $h$--is isomorphic to the blow-up $Bl$ of $\Pr^1 \times \Pr^1$ at the points $\{(0,1) \times (1,0)\}$ and $\{(1,0) \times (0,1)\}$.  Finally, the scheme $\tilde T^{(5)}$ is the disjoint union of six irreducible components (points). They are: $T_{12789}$, $T_{16789}$, $T_{24789}$, $T_{34789}$, $T_{35789}$, $T_{56789}$. Their support maps isomorphically onto the (point) diagonal $\Delta_{123}: Y_{123} \to Y_{123} \times Y_{123}$. 
\end{lem}
\begin{proof} The local description of $X \times_S X$ at $(P,P)$ is given by the equations $z_1z_2z_3 = w_1w_2w_3$ and $z_1z_2z_3 = t$, for $t \in S$ a fixed parameter on the disk. We choose the standard orientation of the sets $\{z_1,z_2,z_3\}$ and $\{w_1,w_2,w_3\}$ and we write $w_i' = \frac{w_i}{z_i}$, $z_i' = \frac{z_i}{w_i}$ for $i = 1,2,3$. After three blows-up of $X \times_S X$ along the subvarieties $z_i = 0 = w_i$, the resulting model $Z$ is non singular as one can see by looking at the first of the following tables which describes $Z$ on each of the eight charts $\mathcal U_j$ who cover it. In the second table, we have collected for each $\mathcal U_j$, the description of the non empty divisors $T_k \in T^{(1)}$ there and the third table shows the ``diagonal'' curves $\delta$ and $\tilde{\delta}$ defined in each chart. The remaining charts describe the pullbacks $p_1^*(\frac{dz_i}{z_i} \wedge \frac{dz_j}{z_j})$, $p_2^*(\frac{dw_i}{w_i} \wedge \frac{dw_j}{w_j})$, $p_1^*(\frac{dz_1}{z_1} \wedge \frac{dz_2}{z_2} \wedge \frac{dz_3}{z_3})$ and $p_2^*(\frac{dw_1}{w_1} \wedge \frac{dw_2}{w_2} \wedge \frac{dw_3}{w_3})$ in terms of the related descriptions by cocycles classes in the corresponding cohomologies.
\begin{center}
\begin{tabular}{c|c} 
Open sets & Loc. coordinates and relations \\ \hline\hline
  {\rule[-3mm]{0mm}{8mm}} 
$\mathcal U_1$ & $\{w_1',w_2',w_3',z_1,z_2,z_3\},~w_1'w_2'w_3' = 1$\\ \hline
  {\rule[-3mm]{0mm}{8mm}} 
$\mathcal U_2$ &$\{w_1',w_2',z_1,z_2,w_3\},~w_1'w_2' = z_3'$  \\ \hline
  {\rule[-3mm]{0mm}{8mm}} 
$\mathcal U_3$ &$\{w_1',w_3',z_1,z_3,w_2\},~w_1'w_3' = z_2'$
  \\ \hline
  {\rule[-3mm]{0mm}{8mm}} 
$\mathcal U_4$ &$\{z_2',z_3',z_1,w_2,w_3\},~z_2'z_3' = w_1'$  \\ \hline
  {\rule[-3mm]{0mm}{8mm}} 
$\mathcal U_5$ &$\{w_2',w_3',z_2,z_3,w_1\},~w_2'w_3' = z_1'$  \\ \hline
  {\rule[-3mm]{0mm}{8mm}} 
$\mathcal U_6$ &$\{z_1',z_3',z_2,w_1,w_3\},~z_1'z_3' = w_2'$  \\ \hline
  {\rule[-3mm]{0mm}{8mm}} 
$\mathcal U_7$ &$\{z_1',z_2',z_3,w_1,w_2\},~z_1'z_2' = w_3'$
  \\ \hline
  {\rule[-3mm]{0mm}{8mm}} 
$\mathcal U_8$ &$\{z_1',z_2',z_3',w_1,w_2,w_3\},~z_1'z_2'z_3' = 1$
  \\ \hline
\end{tabular}
\end{center}\vspace{.1in}

\begin{center}
\begin{tabular}{c|p{9.9cm}} 
Open sets & 
\qquad\qquad\qquad\qquad Divisors\\ \hline\hline
  {\rule[-3mm]{0mm}{8mm}} 
$\mathcal U_1$ & $T_7 = \{z_1 = 0\},~T_8 = \{z_2 = 0\},~T_9 = \{z_3 = 0\}$   \\ \hline
  {\rule[-3mm]{0mm}{8mm}} 
$\mathcal U_2$ & $T_5 = \{w_1' = 0\},~T_6 = \{w_2' = 0\},~T_7 = \{z_1 = 0\}$,\\&{} \qquad\qquad $T_8 = \{z_2 = 0\},~T_9 = \{w_3 = 0\}$  \\ \hline
  {\rule[-3mm]{0mm}{8mm}} 
$\mathcal U_3$ & $T_3 = \{w_1' = 0\},~T_4 = \{w_3' = 0\},~T_7 = \{z_1 = 0\}$,\\&{}\qquad\qquad $T_8 = \{w_2 = 0\},~T_9 = \{z_3 = 0\}$ \\ \hline
  {\rule[-3mm]{0mm}{8mm}} 
$\mathcal U_4$ & $T_3 = \{z_2' = 0\},~T_5 = \{z_3' = 0\},~T_7 = \{z_1 = 0\}$,\\&{}\qquad\qquad $T_8 = \{w_2 = 0\},~T_9 = \{w_3 = 0\}$ \\ \hline
  {\rule[-3mm]{0mm}{8mm}} 
$\mathcal U_5$ & $T_1 = \{w_2' = 0\},~T_2 = \{w_3' = 0\},~T_7 = \{w_1 = 0\}$,\\&{}\qquad\qquad $T_8 = \{z_2 = 0\},~T_9 = \{z_3 = 0\}$ \\ \hline
  {\rule[-3mm]{0mm}{8mm}} 
$\mathcal U_6$ & $T_1 = \{z_1' = 0\},~T_6 = \{z_3' = 0\},~T_7 = \{w_1 = 0\}$,\\&{}\qquad\qquad $T_8 = \{z_2 = 0\},~T_9 = \{w_3 = 0\}$  \\ \hline
  {\rule[-3mm]{0mm}{8mm}} 
$\mathcal U_7$ & $T_2 = \{z_1' = 0\},~T_4 = \{z_2' = 0\},~T_7 = \{w_1 = 0\}$,\\&{}\qquad\qquad $T_8 = \{w_2 = 0\},~T_9 = \{z_3 = 0\}$  \\ \hline
  {\rule[-3mm]{0mm}{8mm}} 
$\mathcal U_8$ & $T_7 = \{w_1 = 0\},~T_8 = \{w_2 = 0\},~T_9 = \{w_3 = 0\}$  \\ \hline
\end{tabular}
\end{center}\vspace{.1in}

\begin{center}
\begin{tabular}{c|p{9.9cm}} 
Open sets & 
\qquad\qquad\qquad ``Diagonal'' curves\\ \hline\hline
  {\rule[-3mm]{0mm}{8mm}} 
$\mathcal U_1$ &\qquad\qquad\qquad\quad none   \\ \hline
  {\rule[-3mm]{0mm}{8mm}} 
$\mathcal U_2$ & $\delta_{13} = \{w_1' = z_1 = w_3 = 0,~w_2' = 1\} \subset T_{579},~\delta_{13} \cap T_8 \neq \emptyset$\\&{}$\delta_{23} = \{w_2' = z_2 = w_3 = 0,~ w_1' = 1\} \subset T_{689},~\delta_{23} \cap T_7 \neq \emptyset$  \\ \hline
  {\rule[-3mm]{0mm}{8mm}} 
$\mathcal U_3$ & $\delta_{12} = \{w_1' = z_1 = w_2 = 0,~w_3' = 1\} \subset T_{378},~\delta_{12} \cap T_{9} \neq \emptyset$\\&{}$\tilde{\delta}_{23} = \{w_3' = z_3 = w_2 = 0,~w_1' = 1\} \subset T_{489},~\tilde{\delta}_{23} \cap T_{7} \neq \emptyset$ \\ \hline
  {\rule[-3mm]{0mm}{8mm}} 
$\mathcal U_4$ & $\delta_{12} = \{z_2' = z_1 = w_2 = 0,~z_3' = 1\} \subset T_{378},~\delta_{12} \cap T_9 \neq \emptyset$\\&{} $\delta_{13} = \{z_3' = z_1 = w_3 = 0,~z_2' = 1\} \subset T_{579},~\delta_{13} \cap T_8 \neq \emptyset$ \\ \hline
  {\rule[-3mm]{0mm}{8mm}} 
$\mathcal U_5$ & $\tilde{\delta}_{12} = \{w_2' = z_2 = w_1 = 0,~w_3' = 1\} \subset T_{178},~\tilde{\delta}_{12} \cap T_9 \neq \emptyset$\\&{} $\tilde{\delta}_{13} = \{w_3' = z_3 = w_1 = 0,~w_2' = 1\} \subset T_{279}
,~\tilde{\delta}_{13} \cap T_8 \neq \emptyset$ \\ \hline
  {\rule[-3mm]{0mm}{8mm}} 
$\mathcal U_6$ & $\tilde{\delta}_{12} = \{z_1' = z_2 = w_1 = 0,~z_3' = 1\} \subset T_{178},~\tilde{\delta}_{12} \cap T_9 \neq \emptyset$\\&{} $\delta_{23} = \{z_3' = z_2 = w_3 = 0,~z_1' = 1\} \subset T_{689},~\tilde{\delta}_{23} \cap T_7 \neq \emptyset$  \\ \hline
  {\rule[-3mm]{0mm}{8mm}} 
$\mathcal U_7$ & $\tilde{\delta}_{13} = \{z_1' = z_3 = w_1 = 0,~z_2' = 1\} \subset T_{279},~\tilde{\delta}_{13} \cap T_8 \neq \emptyset$\\&{} $\tilde{\delta}_{23} = \{z_2' = z_3 = w_2 = 0,~z_1' = 1\} \subset T_{489},~\tilde{\delta}_{23} \cap T_7 \neq \emptyset$  \\ \hline
  {\rule[-3mm]{0mm}{8mm}} 
$\mathcal U_8$ &\qquad\qquad\qquad\qquad none  \\ \hline
\end{tabular}
\end{center}\vspace{.1in}

Denote by $v_{Y_{ij}}$ a class in $H^*(Y_{ij},\C)$ and by $v_{T_{lk}}$ a class in $H^*(\tilde T^{(2)},\C)$. Then we have

\begin{center}
\begin{tabular}{c|c|c} 
Open sets & $p_1^*(v_{Y_{12}})$ & $p_2^*(v_{Y_{12}})$\\ \hline\hline
{\rule[-3mm]{0mm}{8mm}}
$\mathcal U_1$ & $v_{T_{78}}$ & $v_{T_{78}}$\\ \hline
{\rule[-3mm]{0mm}{8mm}} 
$\mathcal U_2$ & $v_{T_{78}}$ & $v_{T_{56}}+v_{T_{58}}-v_{T_{67}}+v_{T_{78}}$\\ \hline
{\rule[-3mm]{0mm}{8mm}} 
$\mathcal U_3$ & $-v_{T_{37}}-v_{T_{47}}+v_{T_{78}}$ & $v_{T_{38}}+v_{T_{78}}$\\ \hline
{\rule[-3mm]{0mm}{8mm}}
$\mathcal U_4$ & $-v_{T_{37}}+v_{T_{78}}$ & $v_{T_{78}}+v_{T_{38}}+v_{T_{58}} $\\ \hline
{\rule[-3mm]{0mm}{8mm}} 
$\mathcal U_5$ & $v_{T_{18}}+v_{T_{28}}+v_{T_{78}}$ & $v_{T_{78}}-v_{T_{17}}$\\ \hline
{\rule[-3mm]{0mm}{8mm}}
$\mathcal U_6$ & $v_{T_{18}}+v_{T_{78}}$ & $v_{T_{78}}-v_{T_{17}}-v_{T_{67}}$\\ \hline
{\rule[-3mm]{0mm}{8mm}}
$\mathcal U_7$ & $v_{T_{24}}+v_{T_{28}}-v_{T_{47}}+v_{T_{78}}$ & $v_{T_{78}}$\\ \hline
{\rule[-3mm]{0mm}{8mm}}
$\mathcal U_8$ & $v_{T_{78}}$ & $v_{T_{78}}$\\ \hline
\end{tabular}
\end{center}\vspace{.1in}

Hence, the global description of the pullbacks $p_1^*(v_{Y_{12}})$ and $p_2^*(v_{Y_{12}})$ are
\begin{align*}
p_1^*(v_{Y_{12}}) &= (v_{T_{18}} + v_{T_{28}} - v_{T_{37}} - v_{T_{47}} + v_{T_{78}}) + v_{T_{24}}\\\\
p_2^*(v_{Y_{12}}) &= (- v_{T_{17}} + v_{T_{38}} + v_{T_{58}} - v_{T_{67}} + v_{T_{78}}) + v_{T_{56}}.
\end{align*}\vspace{.1in}

\begin{center}
\begin{tabular}{c|c|c} 
Open sets & $p_1^*(v_{Y_{13}})$ & $p_2^*(v_{Y_{13}})$\\ \hline\hline
{\rule[-3mm]{0mm}{8mm}}
$\mathcal U_1$ & $v_{T_{79}}$ & $v_{T_{79}}$\\ \hline
{\rule[-3mm]{0mm}{8mm}} 
$\mathcal U_2$ & $-v_{T_{57}}-v_{T_{67}}+v_{T_{79}}$ & $v_{T_{79}}+v_{T_{59}}$\\ \hline
{\rule[-3mm]{0mm}{8mm}} 
$\mathcal U_3$ & $v_{T_{79}}$ & $v_{T_{79}}-v_{T_{47}}+v_{T_{39}}+v_{T_{34}}$\\ \hline
{\rule[-3mm]{0mm}{8mm}}
$\mathcal U_4$ & $-v_{T_{57}}+v_{T_{79}}$ & $v_{T_{79}}+v_{T_{39}}+v_{T_{59}} $\\ \hline
{\rule[-3mm]{0mm}{8mm}} 
$\mathcal U_5$ & $v_{T_{19}}+v_{T_{29}}+v_{T_{79}}$ & $v_{T_{79}}-v_{T_{27}}$\\ \hline
{\rule[-3mm]{0mm}{8mm}}
$\mathcal U_6$ & $v_{T_{16}}+v_{T_{19}}-v_{T_{67}}+v_{T_{79}}$ & $v_{T_{79}}$\\ \hline
{\rule[-3mm]{0mm}{8mm}}
$\mathcal U_7$ & $v_{T_{29}}+v_{T_{79}}$ & $v_{T_{79}}-v_{T_{27}}-v_{T_{47}}$\\ \hline
{\rule[-3mm]{0mm}{8mm}}
$\mathcal U_8$ & $v_{T_{79}}$ & $v_{T_{79}}$\\ \hline
\end{tabular}
\end{center}\vspace{.1in}

Hence we have the global descriptions
\begin{align*}
p_1^*(v_{Y_{13}}) &= (v_{T_{19}} + v_{T_{29}} - v_{T_{57}} - v_{T_{67}} + v_{T_{79}}) + v_{T_{16}}\\\\
p_2^*(v_{Y_{13}}) &= (- v_{T_{27}} + v_{T_{39}} - v_{T_{47}} + v_{T_{59}} + v_{T_{79}}) + v_{T_{34}}.
\end{align*}\vspace{.1in}

\begin{center}
\begin{tabular}{c|c|c} 
Open sets & $p_1^*(v_{Y_{23}})$ & $p_2^*(v_{Y_{23}})$\\ \hline\hline
{\rule[-3mm]{0mm}{8mm}}
$\mathcal U_1$ & $v_{T_{89}}$ & $v_{T_{89}}$\\ \hline
{\rule[-3mm]{0mm}{8mm}} 
$\mathcal U_2$ & $-v_{T_{58}}-v_{T_{68}}+v_{T_{89}}$ & $v_{T_{89}}+v_{T_{69}}$\\ \hline
{\rule[-3mm]{0mm}{8mm}} 
$\mathcal U_3$ & $v_{T_{39}}+v_{T_{49}}+v_{T_{89}}$ & $v_{T_{89}}-v_{T_{48}}$\\ \hline
{\rule[-3mm]{0mm}{8mm}}
$\mathcal U_4$ & $v_{T_{35}}+v_{T_{39}}-v_{T_{58}}+v_{T_{89}}$ & $v_{T_{89}}$\\ \hline
{\rule[-3mm]{0mm}{8mm}} 
$\mathcal U_5$ & $v_{T_{89}}$ & $v_{T_{89}}-v_{T_{28}}+v_{T_{19}}+v_{T_{12}}$\\ \hline
{\rule[-3mm]{0mm}{8mm}}
$\mathcal U_6$ & $-v_{T_{68}}+v_{T_{89}}$ & $v_{T_{89}}+v_{T_{19}}+v_{T_{69}}$\\ \hline
{\rule[-3mm]{0mm}{8mm}}
$\mathcal U_7$ & $v_{T_{49}}+v_{T_{89}}$ & $v_{T_{89}}-v_{T_{28}}-v_{T_{48}}$\\ \hline
{\rule[-3mm]{0mm}{8mm}}
$\mathcal U_8$ & $v_{T_{89}}$ & $v_{T_{89}}$\\ \hline
\end{tabular}
\end{center}\vspace{.1in}

Finally we have
\begin{align*}
p_1^*(v_{Y_{23}}) &= (v_{T_{39}} + v_{T_{49}} - v_{T_{58}} - v_{T_{68}} + v_{T_{89}}) + v_{T_{35}}\\\\
p_2^*(v_{Y_{23}}) &= (v_{T_{19}} - v_{T_{28}} - v_{T_{48}} + v_{T_{69}} + v_{T_{89}}) + v_{T_{12}}.
\end{align*}\vspace{.1in}

Using the above tables we deduce the following

\begin{center}
\begin{tabular}{c|c|c} 
Open sets & $p_1^*(1_{Y_{123}})$ & $p_2^*(1_{Y_{123}})$\\ \hline\hline
{\rule[-3mm]{0mm}{8mm}}
$\mathcal U_1$ & $1_{T_{789}}$ & $1_{T_{789}}$\\ \hline
{\rule[-3mm]{0mm}{8mm}} 
$\mathcal U_1$ & $1_{T_{578}} + 1_{T_{678}} + 1_{T_{789}}$ & $1_{T_{569}} + 1_{T_{589}} - 1_{T_{679}} + 1_{T_{789}}$\\ \hline
{\rule[-3mm]{0mm}{8mm}} 
$\mathcal U_3$ & $- 1_{T_{379}} - 1_{T_{479}} + 1_{T_{789}}$ & $- 1_{T_{348}} + 1_{T_{389}} + 1_{T_{478}} + 1_{T_{789}}$\\ \hline
{\rule[-3mm]{0mm}{8mm}}
$\mathcal U_4$ & $1_{T_{357}} - 1_{T_{379}} + 1_{T_{578}} + v_{T_{789}}$ & $1_{T_{389}} + 1_{T_{589}} + 1_{T_{789}}$\\ \hline
{\rule[-3mm]{0mm}{8mm}} 
$\mathcal U_5$ & $1_{T_{189}} + 1_{T_{289}} + 1_{T_{789}}$ & $1_{T_{127}} -1_{T_{179}} + 1_{T_{278}} + v_{T_{789}}$\\ \hline
{\rule[-3mm]{0mm}{8mm}}
$\mathcal U_6$ & $- 1_{T_{168}} + 1_{T_{189}} + 1_{T_{678}} + 1_{T_{789}}$ & $- 1_{T_{179}} - 1_{T_{679}} + v_{T_{789}}$\\ \hline
{\rule[-3mm]{0mm}{8mm}}
$\mathcal U_7$ & $1_{T_{249}} + 1_{T_{289}} - 1_{T_{479}} + 1_{T_{789}}$ & $1_{T_{278}} + 1_{T_{478}} - v_{T_{789}}$\\ \hline
{\rule[-3mm]{0mm}{8mm}}
$\mathcal U_8$ & $1_{T_{789}}$ & $1_{T_{789}}$\\ \hline
\end{tabular}
\end{center}\vspace{.1in}

We then obtain
\begin{align*}
p_1^*(1_{Y_{123}}) &= (1_{T_{189}} + 1_{T_{249}} + 1_{T_{289}} - 1_{T_{379}} - 1_{T_{479}} + 1_{T_{789}}) - 1_{T_{168}} + 1_{T_{357}} + 1_{T_{578}} + 1_{T_{678}}\\\\
p_2^*(1_{Y_{123}}) &= (- 1_{T_{179}} + 1_{T_{389}} + 1_{T_{569}} + 1_{T_{589}} - 1_{T_{679}} + 1_{T_{789}}) + 1_{T_{127}} + 1_{T_{278}} - 1_{T_{348}} + 1_{T_{478}}.
\end{align*}\vspace{.1in}

Notice that with the exception of $\mathcal U_1$ and $\mathcal U_8$ that are open sets in $\Af^5$ on which only the exceptional components $T_7$, $T_8$ and $T_9$ are non empty, all the remaining charts $\mathcal U_j$ are isomorphic to $\Af^5$ and in each of them one has five components $T_k$ non empty.

On $\mathcal U_3 \cap \mathcal U_4$ the surface $T_{378}$ contains the curve $\delta_{12}$, and on $\mathcal U_5 \cap \mathcal U_6$, $T_{178}$ contains the curve $\tilde{\delta}_{12}$. The curves $\delta_{12}$ and $\tilde{\delta}_{12}$ are different: \ie $T_1 = \emptyset$ on $\mathcal U_3$ and $\mathcal U_4$, but their supports map isomorphically onto the same diagonal $\Delta_{12}: Y_{12} \to Y_{12} \times Y_{12}$.

Similarly, $\mathcal U_2 \cap \mathcal U_4$ contains $\delta_{13}$ whose support maps isomorphically onto $\Delta_{13}$, whereas $\mathcal U_5 \cap \mathcal U_7$ contains $\tilde{\delta}_{13}$, whose support maps still isomorphically onto $\Delta_{13}$: $\delta_{13} \cap \tilde{\delta}_{13} = \emptyset$.

Finally, $\delta_{23} \subset \mathcal U_2 \cap \mathcal U_6$, $\delta_{23} \simeq \Delta_{23}$, while $\tilde{\delta}_{23} \subset \mathcal U_3 \cap \mathcal U_7$, $\tilde{\delta}_{23} \simeq \Delta_{23}$ and $\delta_{23} \cap \tilde{\delta}_{23} = \emptyset$.

The blow-up $Z_1$ of $X \times_S X$ at $z_1 = 0 = w_1$ is the strict transform of $X \times_S X$ in the blow-up of $\Af^6$ along the corresponding linear subvariety. Let $(\tilde z_1, \tilde w_1)$ be a couple of homogeneus coordinates. The exceptional divisor, say $E_1^{(1)}$, is locally a $\Pr^1_{(\tilde z_1,\tilde w_1)}$--bundle over $\{z_1 = 0 = w_1\}$. Then, the intersection  $E_1^{(1)} \cap Z_1$ is locally defined on $E_1^{(1)}$ by $z_2z_3\tilde z_1 - w_2w_3\tilde w_1 = 0$. The blow $E_1^{(2)}$ of $\Pr^1 \times \{z_1 = 0 = w_1\}$ on $\Pr^1 \times \{z_1 = z_2 = w_1 = w_2 = 0\}$ defines the strict transform of $E_1^{(1)}$ after the second blow-up along $\{z_2 = w_2\}$. Said $E_2^{(2)}$ the exceptional divisor of the second blow-up and $(\tilde z_2, \tilde w_2)$  another couple of homogeneus coordinates, one has $E_1^{(2)} \cap E_2^{(2)} = \Pr^1_{(\tilde z_1,\tilde w_1)} \times \Pr^1_{(\tilde z_2,\tilde w_2)} \times \{z_1 = z_2 = w_1 = w_2 = 0\}$. Finally, after the third blowing at $\{z_3 = 0 = w_3\}$ the three exceptional divisors $E_1^{(3)}$, $E_2^{(3)}$ and $E_3^{(3)}$ will intersect the strict transform $Z$ of $X \times_S X$ along the exceptional surface $T_{789}$. This surface is described by the equation $\tilde z_1\tilde z_2\tilde z_3 - \tilde w_1\tilde w_2\tilde w_3 = 0$ in $E_1^{(3)} \cap E_2^{(3)} \cap E_3^{(3)} = \Pr^1_{(\tilde z_1,\tilde w_1)} \times \Pr^1_{(\tilde z_2,\tilde w_2)} \times \Pr^1_{(\tilde z_3,\tilde w_3)} \times \{z_1 = z_2 = z_3 = w_1 = w_2 = w_3 = 0\} = (\Pr^1)^3$, $(\tilde z_3,\tilde w_3)$ being a third couple of homogeneus coordinates. Let consider the projection $T_{789} \to \Pr^1_{(\tilde z_2,\tilde w_2)} \times \Pr^1_{(\tilde z_3,\tilde w_3)}$. The fibre of this map over a given point in the base $(\Pr^1)^2$ is defined by a linear equation as $\alpha z_1 - \beta w_1 = 0$. If either $\alpha$ or $\beta$ (or both) is not zero, then this fibre is reduced to a single point, so the projection map is locally an isomorphism. On the other hand, $\alpha = 0 = \beta$ happens over the two points $(1,0) \times (0,1)$ and $(0,1) \times (1,0)$, where the fibre is a $\Pr^1$. Since $T_{789}$ is non singular, these two copies of $\Pr^1$ are Cartier divisors, so by the universal property of blow-ups the map factors through the blow-up $Bl$ of $(\Pr^1)^2$ at the two points (\ie $T_{789} \to Bl \to (\Pr^1)^2$). It is easy to see from this description that $T_{789} \simeq Bl$. 

It is straighforward to verify from the second table the description of $\tilde T^{(5)}$ on each chart $\mathcal U_j$ and the statement concerning its support.
\end{proof}

The following result generalizes the description of $[N]$ given in theorem~\ref{th0} for double points degenerations.

\begin{thm}\label{thm1} Let $f: X \to S$ be a semistable degeneration of surfaces as we have considered above. With the same notations as in lemma~\ref{lem1}, let $\pi: Bl \to \Pr^1 \times \Pr^1$ be the morphism definying the blow-up of $\Pr^1 \times \Pr^1$ at the points $\{(0,1) \times (1,0)\}$ and $\{(1,0) \times (0,1)\}$, being $Bl \simeq T_{789}$. Let $F_1 = \pi^*(\text{\{pt\}} \times \Pr^1)$ and $F_2 = \pi^*(\Pr^1 \times \text{\{pt\}})$ be the two fundamental fibres and let $E_1$ and $E_2$ be the two exceptional divisors of $\pi$. The following description of $[N] \in \text{Ker}~\rho^{(4)}$ (\cf \eqref{4}) holds:
\[
[N] = a_{178}\tilde{\delta}_{12} + a_{279}\tilde{\delta}_{13} + a_{378}\delta_{12} + a_{489}\tilde{\delta}_{23} + a_{579}\delta_{13} + a_{689}\delta_{23} + \Gamma.
\]

The 1-cycle $\Gamma \subset Bl$ and the (rational) numbers a's are subject to the following requirements:
\[
\Gamma = xF_1 + yF_2 + zE_1 + wE_2, \qquad\text{with}\quad w = z - 1,\quad x,y,z,w \in \Q
\]
\[
a_{178} - a_{378} = a_{279} - a_{579} = a_{489} - a_{689} = 1
\]

\noindent and the relations among them are given by the following set of equalities
\begin{gather*}\begin{split}
&a_{178} = - w,\qquad a_{279} = -(y + w),\qquad a_{378} = - z,\\
&a_{489} = x + z,\qquad a_{579} = -(y + z),\qquad a_{689} = x + w.
\end{split}
\end{gather*}

Furthermore, for those degenerations with $N^2 \neq 0$, the class $[N^2] \in E_1^{0,4}(Z) = H^0(\tilde T^{(5)},\Q)$ (\cf \eqref{5}) can be exhibited as:
\begin{equation*}
\begin{split}
[N^2] = &~b_{12789}T_{12789} + b_{16789}T_{16789} + b_{24789}T_{24789} + b_{34789}T_{34789} + b_{35789}T_{35789} + \\&+ b_{56789}T_{56789}.
\end{split}
\end{equation*}

The (rational) numbers b's must satisfy the following equation:
\[
- b_{12789} + b_{16789} - b_{24789} + b_{34789} - b_{35789} - b_{56789} = 1.
\]

Hence, the induced classes of $[N]$ in $gr^L_2H^4(T,\Q)(1)$ and of $[N^2]$ in $gr^L_0H^4(T,\Q)$ (\ie modulo boundary relations via the restriction maps $\rho^{(3)}$ and $\rho^{(5)}$ \cf~\eqref{diff}) determine algebraic cocycles of dimension one and zero respectively.
\end{thm}

\begin{proof} We will determine $[N]$ as a cocycle making the following square commute (\ie this is the one one has to study for a degeneration of K-3 surfaces of the type mentioned above)
\begin{equation*}
\begin{CD}
gr^L_3H^2(\tsX \times \tsX, \Q) @>{[N]~\cdot}>> gr^L_5H^{6}(\tsX \times \tsX, \Q)(1) = (E_2^{1,5})(1) \\
@A{(p_1)^*}AA @VV{(p_2)_*}V \\
E_2^{-1,3} = gr^L_3H^2(\tsX, \Q) @>{N}>> (gr^L_1H^2(\tsX, \Q))(-1) = (E_2^{1,1})(-1)
\end{CD}
\end{equation*}\vspace{.1in}

Note that besides the commutativity of the square, one has to impose another condition on $[N]$ in order for it to represent the operator $N$. That arises from \eqref{4}. Namely, the representative of $N$ in $(E_1^{2,2})(1) = H^2(\tilde T^{(3)},\Q)(1)$ must belong to the kernel of the related restriction map $\rho^{(4)}$. This condition was automatically satisfied for double point degenerations since $T^{(4)} = \emptyset$ always in that case. We will explicitly describe a representative $[N]$ of $N$ in  $(E_1^{2,2})(1)$ that satisfies the commutativity of the following square
\begin{equation}\label{6}
\begin{CD}
H^1(\tilde T^{(2)},\Q)(-1) @>{[N]~\cdot}>> H^{5}(\tilde T^{(2)},\Q)(1) \\
@A{p_1^*}AA @VV{(p_2)_*}V \\
H^1(\tilde Y^{(2)},\Q)(-1) @= H^1(\tilde Y^{(2)},\Q)(-1).
\end{CD}
\end{equation}\vspace{.1in}

With the notations introduced in lemma~\ref{lem1} we first remark that the cocycles $[\delta_{ij}] = (\Delta_{ij})_*(1_{Y_{ij}})$ ($i,j = 1,2,3$, $i \neq j$), $\Delta_{ij}: Y_{ij} \to Y_{ij} \times Y_{ij}$ being the diagonal embedding, evidently satisfy the cohomological equality
\[
(p_2)_*(\Delta_*(1_{Y_{ij}}) \cdot (p_1)^*(v)) = (p_2)_*(\Delta_*\Delta^*p_1^*(v)) = (p_2)_*(\Delta_*(v)) = v
\]

\noindent for $1_{Y_{ij}} \in H^0(Y_{ij},\Q)$ and any element $v \in H^1(\tilde Y^{(2)},\Q)(-1)$. However, since a simple  linear combination as $a_{178}\tilde{\delta}_{12} + a_{279}\tilde{\delta}_{13} + a_{378}\delta_{12} + a_{489}\tilde{\delta}_{23} + a_{579}\delta_{13} + a_{689}\delta_{23}$ (the coefficients a's are integers) does not satisfy the requirement of being in the kernel of the restriction map $\rho^{(4)}$ (\cf~\eqref{4} and \eqref{diff}), we have to add to the above ``diagonal'' definition a 1-cocycle $\Gamma \subset T_{789}$,  so that the completed linear combination defines an element in $(E_2^{2,2})(1)$ representing $N$. Notice that since the exceptional surface $T_{789}$ projects down via $p_2$, onto the triple point $P$, this modification by $\Gamma$ does not spoil the commutativity of \eqref{6}, once we have checked it for the partial representative of $[N]$ given in terms of the above diagonals.

The 1-cycle $\Gamma$ will be described as a combination of the generators $F_1,~F_2,~E_1,~E_2$ of the Neron-Severi group $NS(T_{789})$. First of all, let consider the six curves $T_{k789}$ for $k = 1,\dots,6$. They are elements of $\tilde T^{(4)}$. We describe them using the generators of $NS(T_{789})$. Because $\pi(T_{1789}) = \{(0,1) \times (1,0)\}$, $T_{1789} = E_2$. Similarly, we have $T_{3789} = E_1$, as $\pi(T_{3789}) = \{(1,0) \times (0,1)\}$. The remaining four curves are described using the projection formula. For example, we know that $\pi(T_{2789}) = (0,1) \times \Pr^1$ and that $\pi^*((0,1) \times \Pr^1) = F_1 = E_2 + T_{2789}$. Hence we have $T_{2789} = F_1 - E_2$. With a similar procedure we obtain $T_{4789} = F_2 - E_1$, $T_{5789} = F_1 - E_1$ and $T_{6789} = F_2 - E_2$. The geometry of the intersections among the generators of $NS(T_{789})$ is well known, namely $E_1 \cdot E_2 = E_1 \cdot F_2 = E_1 \cdot F_1 = E_2 \cdot F_1 = E_2 \cdot F_2 = F_1 \cdot F_1 = F_2 \cdot F_2 = 0$, $E_1 \cdot E_1 = -1 = E_2 \cdot E_2$ and $F_1 \cdot F_2 = 1$.

Let $\Gamma = xF_1 + yF_2 + zE_1 + wE_2$ be an element of $NS(T_{789})$, with $x,y,z,w \in \Q$. Then, we must solve 
\[
[N] = a_{178}\tilde{\delta}_{12} + a_{279}\tilde{\delta}_{13} + a_{378}\delta_{12} + a_{489}\tilde{\delta}_{23} + a_{579}\delta_{13} + a_{689}\delta_{23} + \Gamma
\]

\noindent for $\Gamma$ subject to the condition that $[N]$ is in $\text{ker}~\rho^{(4)}$, for $\rho^{(4)} = \sum_{u=1}^4(-1)^{u-1}\rho_u^{(4)}$ (\cf~\eqref{diff}). For example we have $\rho^{(4)}(a_{178}\tilde{\delta}_{12}) = -a_{178}(\tilde{\delta}_{12} \cdot T_9)$, while $\rho^{(4)}(a_{279}\tilde{\delta}_{13}) = a_{279}(\tilde{\delta}_{13} \cdot T_8)$. Following these rules we obtain the system
\begin{gather}\begin{split}\label{7}
&a_{178} = \Gamma \cdot T_{1789} = - w,\qquad a_{279} = -\Gamma \cdot T_{2789} = -(y + w)\\
&a_{378} = \Gamma \cdot T_{3789} = - z,\qquad a_{489} = \Gamma \cdot T_{4789} = x + z\\
&a_{579} = -\Gamma \cdot T_{5789} = -(y + z),\qquad a_{689} = \Gamma \cdot T_{6789} = x + w.
\end{split}
\end{gather}

For the standard choice of the orientations of $\{z_1,z_2,z_3\}$ and $\{w_1,w_2,w_3\}$ and the numbering of the $T_i$'s setted in lemma~\ref{lem1}, the local description of the pullbacks $\frac{dz_i}{z_i} \wedge \frac{dz_j}{z_j}$ and $\frac{dw_i}{w_i} \wedge \frac{dw_j}{w_j}$ ($i \neq j$, $i,j = 1,2,3$) in terms of cohomology classes $v_{T_{ij}}$ and $v_{T_{ijk}}$, is given following the tables shown in the proof of lemma~\ref{lem1}.

Let $v_{ij} \in H^1(\tilde Y^{(2)},\Q)(-1)$, then via the multiplicative rule described in the Appendix (\cf the similar calculation done in the proof of theorem~\ref{th0}) we obtain 
$$ [N] \cdot p_1^*(v_{12} + v_{13} + v_{23}) =
$$
$$
= [N] \cdot (v_{T_{18}} + v_{T_{78}} + v_{T_{29}} + v_{T_{79}} + v_{T_{49}} + v_{T_{89}}) = 
$$
$$
= a_{178} g_1(\tilde{\delta}_{12} \cdot v_{T_{18}}) - a_{378} g_7(\delta_{12} \cdot v_{T_{78}}) + a_{279} g_2(\tilde{\delta}_{13} \cdot v_{T_{29}}) - a_{579} g_7(\delta_{13} \cdot v_{T_{79}}) + 
$$
$$
+ a_{489}(g_4(\tilde{\delta}_{23} \cdot v_{T_{49}}) - a_{689}(g_8(\delta_{23} \cdot v_{T_{89}}) = 
$$
$$
= a_{178}v_{78}(1) - a_{378}v_{38}(1) + a_{279}v_{79}(1) - a_{579}v_{59}(1) + a_{489}v_{89}(1) - a_{689}v_{69}(1)
$$ 

\noindent where $g_j$ are the pushforward maps defined in the Appendix. Applying the map $(p_2)_*$ we have
$$
(p_2)_*([N] \cdot p_1^*(v_{12} + v_{13} + v_{23})) = (a_{178} - a_{378})v_{12} + (a_{279} - a_{579})v_{13} + (a_{489} - a_{689})v_{23}.
$$

The commutativity of the diagram \eqref{6} is then equivalent to the requirement
\begin{equation}\label{8}
a_{178} - a_{378} = a_{279} - a_{579} = a_{489} - a_{689} = 1
\end{equation}

The linear system \eqref{7} may be then read as $z - w = 1$. Therefore, any curve $\Gamma = xF_1 + yF_2 + zE_1 + wE_2$ satisfying the condition $z - w = 1$ can be used in the description of $[N] \in (E_1^{2,2})(1)$.

The description of $[N^2]$ is similar. For instance, from  proposition~\ref{prop1} we have 
\[
[N^2] \in  gr^L_0H^4(T,\Q) \simeq \frac{H^0(\tilde T^{(5)},\Q)}{\text{Im}(\rho^{(5)}: H^0(\tilde T^{(4)},\Q) \to H^0(\tilde T^{(5)},\Q))}
\]

Via the procedure described in \eqref{2}, $[N^2]$ is then determined in terms of the commutativity of the following square
\begin{equation*}
\begin{CD}
gr^L_4H^2(\tsX \times \tsX, \Q) @>{[N^2]~\cdot}>> gr^L_4H^{6}(\tsX \times \tsX, \Q) = E_2^{2,4} \\
@A{(p_1)^*}AA @VV{(p_2)_*}V \\
E_2^{-2,4} = gr^L_4H^2(\tsX, \Q) @>{N^2}>> (gr^L_0H^2(\tsX, \Q))(-2) = (E_2^{2,0})(-2).
\end{CD}
\end{equation*}

The related $E_1$ description is
\begin{equation*}
\begin{CD}
H^0(\tilde T^{(3)},\Q)(-2) @>{[N^2]~\cdot}>> H^{4}(\tilde T^{(3)},\Q) \\
@A{p_1^*}AA @VV{(p_2)_*}V \\
H^0(\tilde Y^{(3)},\Q)(-2) @= H^0(\tilde Y^{(3)},\Q)(-2).
\end{CD}
\end{equation*}

The scheme $\tilde T^{(5)}$ is the disjoint union of the zero dimensional schemes $T_{12789}$, $T_{16789}$, $T_{35789}$ and $T_{56789}$. Their support map all isomorphically onto the diagonal $\Delta_{123}: Y_{123} \to Y_{123} \times Y_{123}$. With a similar procedure as the one used above to describe $[N]$, we write
\begin{equation*}
\begin{split}
[N^2] = &~b_{12789}T_{12789} + b_{16789}T_{16789} + b_{24789}T_{24789} + b_{34789}T_{34789} + b_{35789}T_{35789} + \\&+ b_{56789}T_{56789}
\end{split}
\end{equation*}

\noindent for some integers $b$'s. Imposing the commutativity of the above diagram, by means of the description of the pullbacks $p_1^*(1_{Y_{123}})$ and $p_2^*(1_{Y_{123}})$ as shown in the last table appearing in the proof of lemma~\ref{lem1}, we finally get the condition
\[
- b_{12789} + b_{16789} - b_{24789} + b_{34789} - b_{35789} - b_{56789} = 1.
\]
\end{proof}\vspace{.1in}

It is straightforward to verify that both $[N]$ and $[N^2]$ make diagrams like \eqref{2} commute, for any choice of the indices $*$ and $r$.\vspace{.1in}

\begin{rem}\end{rem}

It is easy to verify that the description of $[N]$ and $[N^2]$ given in theorem~\ref{thm1} holds also for a normal--crossings degeneration (not semistable) like $f(z_1,\ldots,z_n) = z_i^2z_j$, $i,j \in [1,n]$, $i \neq j$. This applies in particular to the case of normal--crossings degenerations of curves with triple points as described above. The desingularization process of the threefold $X \times_S X$ is obtained via two blow-ups along $z_i = 0 = w_i$ and $z_j = 0 = w_j$ by analogy to what we have done in Remark~\ref{rem1}. For the description of $[N]$ we also refers to the same Remark.\vspace{.2in}

\section{An arithmetic interpretation of the monodromy operator in mixed characteristic}\label{f}

The calculations on the  geometric description of $[N^i]$ that we have done in the previous sections only involve the (local) geometry of the special fibre of a degeneration. Hence they equally hold in mixed characteristic also, \ie for a degeneration $f: \mathcal X \to \text{Spec}(\Lambda) = S$, where $\Lambda$ is a Henselian discrete valuation ring with $\eta$ and $v$ as its generic and closed points respectively. In analogy with the classical case, the model $\mathcal X$ is assumed to be proper and non singular and the map $f$ is supposed to be flat, smooth over the generic point $\eta$ and with a normal--crossings special fibre $Y$ defined over the {\it finite} field $k(v)$ of characteristic $p > 0$.

Locally, for the \'etale topology $\mathcal X$ is $S$-isomorphic to $S[x_1,\ldots,x_n]/(x_1^{e_1}\cdots x_k^{e_k} - \pi)$, where $\pi$ is a uniformizing parameter in $\Lambda$ and $e_i \in \Z,~\forall i = 1,\ldots,k$. For simplicity, we also assume that $\Lambda$ is a {\it finite extension} of $\Zl$ or $\Ql$, where $l \neq p$ is a prime number.

The complex of nearby cycles is then defined as $\R\P(\Lambda) := \bar{i}^{-1}\R\bar{j}_*\Lambda$. Here $i: Y \to \mathcal X$ (\resp $j: \mathcal X_{\eta} \to \mathcal X$) is the natural closed (\resp open) embedding that one ``extends'' to the algebraic closure $k(\bar v)$ of $k(v)$ (\resp a separable closure $k(\bar\eta)$ of $k(\eta)$). Assume that the multiplicities $e_i$ are prime to $\ell$ and $g.c.d.(e_i,p) = 1$. Then, the wild inertia acts trivially on $\R\P(\Lambda)$ and the theory exposed in \cite{RZ} shows that the nearby cycle complex has an abstract description in the derived category $D^+(Y,\Lambda [\Zl(1)])$ of the abelian category of complexes of sheaves of $\Lambda [\Zl(1)]$--modules on $Y$, by a complex $A^\bu_{\mathcal X,\Lambda}$, supported on $Y$. $A^\bu_{\mathcal X,\Lambda}$ can be interpreted as the analogue of the Steenbrink resolution in the classical case. Therefore, the related study of it goes in parallel with the classical one in equal characteristic zero. We refer to \op and \cite{I} (\eg Th\'eor\`eme~3.2) for further detail. \vspace{.1in}

The power maps ($n \in [0,2d]$, $i \ge 0$, $d = \dim~\mathcal X_{\eta}$)
\[
N^i: H^n(\mathcal X_{\bar\eta},\Lambda) \to H^n(\mathcal X_{\bar\eta},\Lambda)(-i)
\]

\noindent define elements
\[
N^i \in \bigoplus_{n\ge 0}\biggl [H^{2d-n}(\mathcal X_{\bar\eta},\Lambda)(d) \otimes H^n(\mathcal X_{\bar\eta},\Lambda)(-i)\biggr ]^G = \biggl [H^{2d}(\mathcal X_{\bar\eta} \times \mathcal X_{\bar\eta},\Lambda)(d-i)\biggr ]^G
\]

\noindent invariant for the action of the Galois group $G = \text{Gal}(\bar\eta/\eta)$ on the cohomology of the product $\mathcal X_{\bar\eta} \times \mathcal X_{\bar\eta}$.  Assume that $f: \mathcal X \to S$ has at worst triple points. Then, the singularities of both $\mathcal X \times_S \mathcal X$ and $Y \times Y$ can be resolved locally around each singular point by a sequence of at most three blows-up, as we described in details in $\S\S$~\ref{a},\ref{c},\ref{e}. The resulting degeneration $h: \mathcal Z \to S$ is normal--crossings with special fibre $T = T_1 \cup\ldots\cup T_N$. Let $\mathcal X_{\bar\eta} \times \mathcal X_{\bar\eta} = \mathcal Z_{\bar\eta}$ be its geometric generic fibre. Denote by $\tilde N = 1 \otimes N + N \otimes 1$ the logarithm of the local monodromy on the product degeneration $h$. Then, the analogue of proposition~\ref{prop1} is the following
\begin{prop}\label{prop2} Assume the monodromy-weight conjecture on $H^*(\mathcal Z_{\bar\eta},\Lambda)$. Then
\[
N^i \in \text{Ker}(\tilde N) \cap H^{2d}(\mathcal Z_{\bar\eta},\Lambda(d-i))^{F=1} \simeq
\]
\[
\simeq \text{Ker}(\tilde N) \cap \biggl ((gr^L_{2(d-i)}H^{2d}(\mathcal Z_{\bar\eta},\Lambda))(d-i)\biggr )^{F=1} \simeq \biggl ((gr^L_{2(d-i)}H^{2d}(T,\Lambda))(d-i)\biggr )^{F=1} \simeq
\]
\[
\simeq \biggl [\frac{\text{Ker}(\rho^{(2(i+1)}: H^{2(d-i)}(\tilde{T}^{(2i+1)},\Lambda)(d-i) \to H^{2(d-i)}(\tilde{T}^{(2(i+1))},\Lambda)(d-i))}{\text{Image}~\rho}\biggr ]^{F=1}
\]

\noindent where $F$ is the geometric Frobenius.
\end{prop}

The following result shows the relation of proposition~\ref{prop2} with the arithmetic of the degeneration $h$
\begin{thm} Assume the monodromy-weight conjecture on $\mathcal Z_{\bar\eta}$ and the semisimplicity of the action of the frobenius $F$ on $H^*(\mathcal Z_{\bar\eta},\Lambda)^I$. Then, for $i > 0$ and $d = \dim~\mathcal X_{\bar\eta}$
$$
\operatorname*{\text{ord}}_{s = d-i}\det(Id - FN(v)^{-s} | H^{2d}(\mathcal Z_{\bar\eta}, \Lambda)^{I}) = 
$$
\[
\text{rk}\biggl [\frac{\text{Ker}(\rho^{(2(i+1)}: H^{2(d-i)}(\tilde{T}^{(2i+1)},\Lambda)(d-i) \to H^{2(d-i)}(\tilde{T}^{(2(i+1))},\Lambda)(d-i))}{\text{Image}~\rho}\biggr ]^{F=1}.
\]

$N(v)$ is the number of elements of the finite residue field $k(v)$.
\label{th1}
\end{thm}
\begin{proof} \cf~\cite{C}, theorem 3.5.
\end{proof}

This result explains geometrically the pole of the local factor at $v$ of the L-function $L(H^{2d}(\mathcal Z_{\bar\eta},\Ql), s)$ at the points $s = d - 1$ and $s = d - 2$, with the presence of the ``diagonal'' cycles representing the monodromy powers on the strata of $T$ as we previously described.\vspace{.2in}

\section{Appendix (by Spencer Bloch)}

Our objective in this appendix is to define a multiplication between the total complex of
$E_1$-terms of the Steenbrink spectral sequence and the graded complex 
\begin{equation}H^*(Y^{(\bullet)}),\ \rho = \text{restriction}
\end{equation}
which is the $E_1$ complex converging to the cohomology of the special
fibre $Y$. We order the components $Y=Y_1\cup\ldots\cup Y_N$ and write $a_{i_0,\dotsc,i_m}\in
H^*(Y_{i_0,\dotsc,i_m},\Q)$.  
The $E_1$-terms of the Steenbrink spectral sequence can be arrayed in a triangular diagram
(compare
\cite{I}, (2.3.8.1)) where each $\bullet$ denotes some $H^*(Y^{(m)},\Q(n))$. 
\begin{equation}\label{3a}\begin{array}{cccccc}\bu &&&&& \\
\uparrow &\nwarrow &&&& \\
\bu & \leftarrow &\bu &&& \\
\uparrow &\nwarrow & \uparrow & \nwarrow && \\
\bu & \leftarrow &\bu & \leftarrow &\bu & \\
\vdots &\makebox[0cm][c]{\,\, $\diagup_{{}^{\text{wt}0}}$} & \vdots &\makebox[0cm][c]{\,\,
$\diagup_{{}^{\text{wt}1}}$}& \vdots & \ddots \\
\bu & \leftarrow & \bu & \leftarrow & \bu & \leftarrow\bu
\end{array}
\end{equation}
Here the horizontal arrows are Gysin maps and the vertical arrows are restriction maps. The
diagonal arrows are (upto twist) the maps $N$ which, on the level of $E_1$ are either the
identity or $0$. The Steenbrink $E_1$-terms, i.e.
the $H^*(Y,\text{gr}_r^L\R\P(\Q))$, are direct sums of terms on a NE-SW diagonal, with weight
$r$ meeting the ''$x$-axis'' at $x=r$. The complex $H^*(Y^{(\bullet)})$ is embedded as the left
hand column, and the resulting multiplication on it is the usual (associative) product 
\begin{equation} a_{i_0,\dotsc,i_m}\otimes b_{j_0,\dotsc,j_n} \mapsto \begin{cases} 0
& i_m \neq j_0 \\ (a\cdot b)_{i_0,\dotsc,i_m,j_1,\dotsc,n_n} & i_m=j_0 \end{cases}
\end{equation} 
The bottom row is a quotient complex calculating the homology of the closed fibre $H_*(Y)$ (with
appropriate twist).  Our
multiplication induces an action of the left hand column on the bottom row, which we
will show induces the cap product (\cite{Sp}, p. 254)
\begin{equation} H^q(Y)\otimes H_n(Y) \to H_{n-q}(Y).
\end{equation}
This module structure, unifying and extending the classical cocycle calculations for cup and cap
product, is of independent interest. Quite possibly it can be extended to a product on the
whole $E_1$-complex, but the daunting sign calculations involved have prevented us from working it
out. 

We will apply this construction to calculate the product 
\begin{equation}[N^i]~\cdot : \bold H^*(T, gr_r^L A_{Z,\Q}^\bu) 
\to \bold H^{*+2d}(T, gr_{r-2i}^L A_{Z,\Q}^\bu(d-i)) 
\end{equation}
from \eqref{2}. 

We return to the situation in section~\ref{a}. In particular, $Z\to X\times_S X$ is a resolution, and
$T\subset Z$ is the special fibre, which we assume is a normal crossings divisor. We write
$E_1(Z)$ for the Steenbrink spectral sequence associated to the degeneration $Z/S$. 
\begin{lem} There exists a class $[N^i]$ in $E_1(Z)$ satisfying
\begin{enumerate}
\item $d_1[N^i]=0$, and the induced class in $E_2$ is the $i$-th power of the monodromy operator
$$N^i \in \text{gr}^L_{2(d-i)}H^{2d}(\tilde X^* \times \tilde X^*,\Q(d-i))
$$
\item $N([N^i])=0$, i.e. in the diagram \eqref{3a}, $[N^i]$ lies in the left hand vertical
column. 
\end{enumerate}
\end{lem}
\begin{proof}We see from proposition~\eqref{2} that the class of $N^i$ is killed by $N$ in $E_2(Z)$.
Let $M$ denote the map on $E_1$ which is inverse to $N$ insofar as possible, i.e. $M$ maps down
and to the right in diagram \eqref{3a}. $M$ is zero on the bottom line. Let $x\in E_1$ represent
$N^i$ in
$E_2$. Then
$Nx = d_1y$. (Here $d_1=d'+d''$ is the total differential.) Since $N$ commutes with $d'$ and
$d''$, and $Nx$ has no term on the bottom row, it follows that $[N^i] := x-d_1My$ is supported on
the left hand column, i.e. killed by $N$. 
\end{proof}

Here is some notation. The special fibre will be $Y=\bigcup
Y_i$, with
$0\le i\le N$. Write
$H^*(Y)$ for cohomology in some fixed constant ring like $\Z$ or $\C$. 
$$I = \{i_0,\ldots,i_m \};\quad J=\{j_0,\ldots,j_n \}\quad\text{(strictly ordered)};\quad Y_I =
\bigcap_{i_k\in I}Y_{i_k}
$$
We will say the pair $I,J$ is {\it admissible} if
$$\exists p \text{ such that } i_m = \max(I) = j_p \text{ and } \{j_0,\ldots,j_p\}\subset I.
$$
In this case, write $j_0=i_{b_0},\ldots, j_{p-1}=i_{b_{p-1}}$. Define
$$a(I,J) := b_0+\ldots+b_{p-1}+mp.
$$
With $I,J$ admissible as above, write
$$J' = \{j_0,\ldots,j_p\};\  J'' = \{j_p,\ldots,j_n\};\  J=J'\cup J'';
\  J'\cap J'' = \{j_p\}=\{i_m\} 
$$
Write
$$I' = J';\ I'' =(I-J')\cup\{i_m\};\ I=I'\cup I'';\ \{i_m\} = I'\cap I'' 
$$
Let $K=I''\cup J''$, and define
\begin{gather}\label{product}\theta(I,J) : H^\alpha(Y_I)\otimes H^\beta(Y_J) \to H^{\alpha+\beta+2p}(Y_K) \\
\theta(I,J)(x\otimes y) := (-1)^{a(I,J)}g_{j_0}\circ\cdots g_{j_{p-1}}(x\cdot y).
\end{gather}
Here $x\cdot y\in H^{\alpha+\beta}(Y_{I\cup J})$, the $g_j$ are Gysin maps, and 
$$g_{j_0}\circ\cdots g_{j_{p-1}}: H^*(Y_{I\cup J}) \to H^{*+2p}(Y_{I''\cup J''}).
$$
If the pair $I,J$ is not admissible, define $\theta(I,J)=0$. 
Define for $I$ as above and $0\le k\le N$
$$\sigma(I,k) := \#\{i\in I\ |\ i<k\}
$$
For $k\notin I$ we have the restriction $\text{rest}_k : H^*(Y_I) \to H^*(Y_{I\cup \{k\}})$.
Define
$$d' := \sum_{k\notin I}(-1)^{\sigma(I,k)}\text{rest}_k : H^*(Y_I) \to \underset{k\notin
I}{\bigoplus}H^*(Y_{I\cup\{k\}})
$$ 
Similarly, for $k\in I$ we have the Gysin $g_k : H^*(Y_I) \to H^{*+2}(Y_{I-\{k\}})$. We define
$$d'' = \sum_{k\in I}(-1)^{\sigma(I,k)}g_k : H^*(Y_I) \to \underset{k\in
I}{\bigoplus}H^{*+2}(Y_{I-\{k\}}).
$$
\begin{thm} With notation as above ($I,J$ not necessarily admissible) the following diagram is
commutative:
$$
\begin{CD}H^*(Y_I)\otimes H^*(Y_J) @>\theta(I,J)>> H^*(Y_K) \\
@VV d'\otimes 1+(-1)^m1\otimes(d'+d'')V  @VVd'+d'' V \\
\underset{\tilde I,\tilde J }{\bigoplus}H^*(Y_{\tilde I})\otimes 
H^*(Y_{\tilde J}) @>\theta(\tilde I,\tilde J)>> \underset{\tilde K = 
\tilde I''\cup \tilde J''}{\bigoplus}H^*(Y_{\tilde K})
\end{CD}
$$
\end{thm}

\begin{rem}\end{rem} A priori the theorem does not suffice to determine the desired mapping
$$ H^*(Y^{\bu})\otimes E_1 \to E_1\quad a\otimes b \mapsto a*b
$$
because a given $H^*(Y_K)$ occurs many times in the diagram \eqref{2} (at every point along a NW
pointing diagonal). However, if we add the condition that the weights (SW-NE diagonals in
\eqref{3a}) should be added, the mapping is defined. It has the property that
$$a*Nb = N(a*b)
$$ 
In particular, there is an induced action on $E_1/NE_1$ which we identify with the bottom row in
\eqref{3a}. This simple complex calculates $H_*(Y)$, and the product coincides with the cap
product. To see this, one notes that the product is correct for two elements in weight $0$, and
that if each $H^*(Y_I)$ is replaced by $\Z$, the acyclic model theorem (\cite{Sp}, p.
165) can be applied.\vspace{.1in} 

\begin{proof}[proof of theorem] The proof consists of many separate cases. In each case we will
check the sign carefully (this is the delicate part) and omit checking that the maps coincide
set-theoretically (which is straightforward). \newline

\noindent {\it case:} $i_m\notin J$. \newline

\noindent In this case, the pair $I,J$ is not admissible, so $\theta(I,J)=0$. We must show
\begin{equation}\label{11}\underset{\tilde I,\tilde J}{\oplus}\theta(\tilde I,\tilde J)\circ
(d'\otimes 1+(-1)^m1\otimes(d'+d''))=0.
\end{equation}
We may ignore non-admissible $\tilde I, \tilde J$. The only way admissible $\tilde I, \tilde J$
can occur in this situation is if for some $p\ge 0$ we have $j_{p-1}<i_m<j_p$ and
$\{j_0,\ldots,j_{p-1}\}\subset I$. (If a subscript for $j$ doesn't fall in $\{0,\ldots,n\}$,
ignore it, i.e. take $j_{-1}=-\infty,\ j_{n+1}=+\infty$.) Assume these conditions hold. Then the
pair $I\cup \{j_p\},\ J$ is admissible and occurs in the image of $d'\otimes 1$. Also the pair
$I,\ J\cup\{i_m\}$ is admissible and occurs in the image of $(-1)^m(1\otimes d')$. We must show
these two contributions cancel. Suppose $j_0=i_{b_0},\ldots,j_{p-1}=i_{b_{p-1}}$. Then the sign
condition we need to verify is 
\begin{multline*}\sigma(I,j_p)+b_0+\cdots+b_{p-1}+p(m+1)\equiv \\
1+m+\sigma(J,i_m)+b_0+\cdots+b_{p-1}+pm
\mod(2)
\end{multline*}
This is correct because $\sigma(I,j_p)=m+1$ and $\sigma(J,i_m)=p$.\newline

\noindent {\it case:} $i_m=j_p\in J,\ \{j_0,\ldots,j_{p-1}\}\not\subset J$. \newline

\noindent This is the other case where $I,J$ is not admissible, so $\theta(I,J)=0$. To get
admissible $\tilde I, \tilde J$ we must have 
$$\exists k,\ 0\le k\le p-1\ \text{such that } j_k\notin I,\ \{j_0,\ldots,\hat
j_k,\ldots,j_{p-1}\}\subset I.
$$ 
Assume this. Then the pairs $(I\cup \{j_k\},J)$ and $(I,J-\{k_k\})$ are admissible. The first
occurs in $\theta(I\cup\{j_k\},J)\circ (d'\otimes 1)$ and the second in
$(-1)^m\theta(I,J-\{j_k\})\circ 1\otimes d''$. The necessary sign condition for cancellation is
$$\sigma(I,j_k)+a(I\cup\{j_k\},J)\stackrel{?}{\equiv} m+1+k+a(I,J-\{j_k\}) \mod(2).
$$
To check this sign condition write $j_r = i_{b_r}$ for $0\le r\le p-1,\ r\not=k$. Then 
\begin{gather*}a(I,J-\{j_k\}) = b_0+\cdots+b_{k-1}+b_{k+1}+\cdots+b_{p-1} +(p-1)m \\
a(I\cup\{j_k\},J) = b_0+\cdots+b_{k-1}+\sigma(I,j_k)+(b_{k+1}+1)+ \\
+\cdots+(b_{p-1}+1)+p(m+1).
\end{gather*}
This yields the necessary congruence. \newline

\noindent For the rest of the proof we assume $I,J$ is admissible. We examine the various terms
in \eqref{11} and show they occur with the same signs in $(d'+d'')\circ\theta(I,J)$. We first
consider terms coming from $d'\otimes 1$, so the target is labelled by $\tilde I=I\cup\{k\},\
\tilde J=J$. 
\newline

\noindent {\it case:} $k<i_m = j_p$. In this case, since $j_p=\min J''$ and $k\notin I\supset
J'$, we have $k\notin J$. The pair $\tilde I = I\cup \{k\},\ J$ is admissible with $\tilde I'' =
I''\cup \{k\}$ and the same decomposition $J=J'\cup J''$. Let $\tilde K = \tilde I''\cup J''=
K\cup \{k\}$. Since $k<j_p = \min J''$, we have
$$\sigma(K,k)=\sigma(I'',k)=\sigma(I,k)-\sigma(J',k)
$$
What we must show, therefore, is that
$$a(I,J)-a(\tilde I,J) \equiv \sigma(J',k) \mod (2)
$$
Write
\begin{gather*}\tilde I = \{\tilde i_0,\ldots,\tilde i_{m+1}\};\ j_0 = \tilde i_{\tilde
b_0},\ldots, j_{p-1} = \tilde i_{\tilde b_{p-1}};\\ 
a(\tilde I,J) = \tilde b_0+\cdots+\tilde b_{p-1}+(m+1)p \\
I = \{i_0,\ldots,i_m\};\ j_r = i_{b_r}, 0\le r\le p-1 \\
a(I,J) = b_0+\cdots+b_{p-1}+mp
\end{gather*}
where
$$\tilde b_\ell = \begin{cases}b_\ell & i_{b_\ell}<k \\ b_\ell +1 & i_{b_\ell} >k \end{cases}
$$
Thus
\begin{multline*}a(\tilde I,J) - a(I,J) = p-\#\{j\in J'-\{j_p\}\ |\ j>k\} = \\
\#\{j\in J'\ |\ j<k\} = \sigma(J',k).
\end{multline*}
This is the desired congruence.\newline

\noindent We continue to consider the contribution of $d'\otimes 1$ with $I,J$
admissible.\newline

\noindent {\it case:} $k>i_m,\ k\not= j_{p+1}$. \newline

\noindent In this case $I\cup\{k\},J$ is not admissible so $\theta(\tilde I,J)=0$. \newline

\noindent {\it case:} $k=j_{p+1}$.

Here $\tilde I:= I\cup \{k\}, \tilde J:= J$ is admissible with 
\begin{gather*}\tilde J' = \{j_0,\ldots,j_{p+1}\}=J'\cup\{k\}=J'\cup\{j_{p+1}\} \\
\tilde J'' = \{j_{p+1},\ldots,j_n\}=J''-\{j_p\};\ \tilde K = \tilde I''\cup\tilde J''=K-\{j_p\}
\end{gather*}
Note in this case $k>i_m$ so $\sigma(I,k)=m+1$. The claim is here that the diagram
$$\begin{CD}H^*(Y_I)\otimes H^*(Y_J) @>\theta(I,J)>> H^*(Y_K) \\
@VV(-1)^{m+1}\text{rest.}\otimes 1 V @VV(-1)^{\sigma(K,j_p)}\text{Gysin}_{j_p}V \\
H^*(Y_{\tilde I})\otimes H^*(Y_J) @>\theta(\tilde I,J)>> H^*(Y_{\tilde K})
\end{CD}
$$
commutes. Note that the right hand vertical arrow (with the sign) is part of $1\otimes d''$. To
verify the signs we need
$$a(I,J)+\sigma(K,j_p)\equiv m+1+a(\tilde I,J).
$$
Since $K=I''\cup J''$ and $k=\max(I'')=\min(J'')$ it is clear that 
$$\sigma(K,j_p) = \#I''-1=m-p.
$$
Also $j_p=i_m$ so with the usual notation $j_r=i_{b_r}$ we get
$$a(\tilde I,J) = b_0+\cdots+b_{p-1}+m+(m+1)(p+1).
$$
Now the desired congruence becomes
$$b_0+\cdots+b_{p-1}+pm+m-p\equiv b_0+\cdots+b_{p-1}+m+(m+1)(p+1)+m+1
$$
This is correct.\newline

\noindent We now consider terms occurring in $(-1)^m(1\otimes d')$ on the left of the diagram in
the statement of the theorem. We assume given $k\notin J$.\newline

\noindent {\it case:} $k>j_p$. \newline

\noindent Note in this case $k\notin I$. Taking $\tilde J=J\cup\{k\},\ \tilde K = K\cup\{k\}$, I
claim the diagram below is commutative:
$$\begin{CD}H^*(Y_I)\otimes H^*(Y_J) @>\theta(I,J)>> H^*(Y_K) \\
@VV(-1)^{m+\sigma(J,k)}1\otimes\text{rest}V @VV(-1)^{\sigma(K,k)}\text{rest}V \\
H^*(Y_I)\otimes H^*(Y_{\tilde J}) @>\theta(I,\tilde J)>> H^*(K_{\tilde K}) 
\end{CD}
$$
(In other words, the contribution in this case is to $d'$ on the right.) Set
$$\tilde J = J'\cup\tilde J'';\ \tilde J'' = J''\cup\{k\};\ K=I''\cup J''.
$$
We have
\begin{gather*} a(I,J)=a(I,\tilde J) \\
\sigma(J,k) = \sigma(J'',k)+p+1 \\
\sigma(K,k)=\sigma(J'',k)+\#I'' = \sigma(J'',k)+m+1-p
\end{gather*}
It follows that
$$m+\sigma(J,k)+a(I,\tilde J)\equiv \sigma(K,k)+a(I,J) \mod(2)
$$
which is the desired sign relation in this case.\newline

\noindent {\it case:} $k<j_p,\ k\notin I$.\newline

\noindent In this case, the pair $I, J\cup\{k\}$ is not admissible, so the contribution is
zero.\newline

\noindent {\it case:} $k<j_p,\ k\in I$.\newline

In this case the pair $I,\tilde J$ is admissible with
\begin{gather*}\tilde J:= J\cup\{k\}=\tilde J'\cup J'';\ \tilde J' = J'\cup\{k\} \\
I=\tilde I = \tilde J'\cup \tilde I'';\ \tilde I''=I''-\{k\};\ \tilde K=K-\{k\}=\tilde I''\cup
J'' 
\end{gather*}
The term in question contributes to $d''$ on the right, and the diagram which commutes is:
$$ \begin{CD} H^*(Y_I)\otimes H^*(Y_J) @>\theta(I,J)>> H^*(Y_K) \\
@VV(-1)^{m+\sigma(J,k)}\text{rest} V @VV (-1)^{\sigma(K,k)}\text{Gysin}_k V \\
H^*(Y_I)\otimes H^*(Y_{\tilde J}) @>\theta(I,\tilde J) >> H^*(Y_{\tilde K})
\end{CD}
$$
The signs will be correct if
$$a(I,J)+\sigma(K,k) \equiv m+\sigma(J,k)+\theta(I,\tilde J) \mod(2)
$$
Write $\tilde J = \{\tilde j_0,\ldots,\tilde j_{m+1}\}$ and $\tilde j_r = i_{\tilde b_r},\ r\le
p$. The desired congruence reads
$$b_0+\cdots+b_{p-1}+mp+\sigma(K,k) \stackrel{?}{\equiv} m+\sigma(J,k)+\tilde b_0+\cdots+\tilde
b_p+(p+1)m
$$
We have
$$\tilde b_\ell = \begin{cases} b_\ell & \ell < \sigma(J',k) \\
\sigma(I,k) & \ell = \sigma(J',k) \\
b_{\ell -1} & \ell > \sigma(J',k) \end{cases}
$$
The condition becomes
$$\sigma(K,k)\stackrel{?}{\equiv} \sigma(J,k)+\sigma(I,k)=
\sigma(J',k)+\sigma(J',k)+\sigma(I'',k),
$$
which is true.\newline

\noindent Finally we consider terms coming from $(-1)^m(1\otimes d'')$ in the lefthand vertical
arrow in the diagram of the theorem. In what follows $j\in J$. \newline

\noindent {\it case:} $j\in J'',\ j\not= j_p$. Define
$$\tilde J=J-\{j\};\ K=I''\cup J'';\ \tilde K=K-\{j\} = I''\cup \tilde J''.
$$
The diagram which commutes is:
$$\begin{CD}
H^*(Y_I)\otimes H^*(Y_J) @>\theta(I,J) >> H^*(Y_K) \\
@VV 1\otimes (-1)^{m+\sigma(J,j)}\text{Gysin}_j V  @VV (-1)^{\sigma(K,j)}\text{Gysin}_j V \\
H^*(Y_I)\otimes H^*(Y_{\tilde J}) @>\theta(I,\tilde J) >> H^*(Y_{\tilde K})
\end{CD}
$$
The sign condition to be checked is
$$m+\sigma(J,j)+a(I,\tilde J) \stackrel{?}{\equiv} a(I,J)+\sigma(K,j) \mod(2).
$$
Our conditions imply $j>j_p$ so $a(I,J)=a(I,\tilde J)$. Also, 
$$\#I'' + \# J' = m+2 \equiv m\mod(2),
$$
so
\begin{gather*} \sigma(K,j) = \#I'' +\sigma(J'',j)-1 \equiv m+ \# J'+\sigma(J'',j)-1 \\
\sigma(J,j) = \sigma(J',j)+\sigma(J'',j)-1 = \# J' -1+\sigma(J'',j).
\end{gather*}
This is the desired condition. \newline

\noindent {\it case:} $j=j_p$.\newline

\noindent In this case, $I, J-\{j\}$ is not admissible, so we get no contribution. \newline

\noindent {\it case:} $j\in J,\ j< j_p$. \newline

\noindent In this case, $j\in J',\ j\not= j_p$. Set
\begin{gather*}\tilde J = J-\{j\};\ \tilde J' = J'-\{j\}; \tilde J'' = J'' \\
\tilde I = I;\ \tilde I'' = I''\cup \{j\}; \ I=\tilde I = \tilde J'\cup \tilde I'' \\
K=I''\cup J'';\ \tilde K = \tilde I''\cup \tilde J'' = K-\{j\}.
\end{gather*}
The sign condition to show we gat a contribution to $d''$ on the right is
$$a(I,J)+\sigma(K,j)\stackrel{?}{\equiv} m+\sigma(J,j)+a(\tilde I, \tilde J) \mod(2).
$$
Writing $j=j_\ell = i_{b_\ell}$ the condition becomes
\begin{multline*}b_0+\cdots+b_{p-1}+mp+\sigma(K,j)\stackrel{?}{\equiv} \\
b_0+\cdots+\hat b_\ell
+b_{\ell +1}+\cdots +b_{p-1}+m(p-1)+m+\sigma (J,j)
\end{multline*}
This is true because
\begin{gather*}b_\ell = \sigma(I,j) = \sigma(I'',j)+\sigma(J',j) \\
\sigma(K,j)=\sigma(I'',j);\quad \sigma(J,j)=\sigma(J',j).
\end{gather*}

The proof is completed by checking that all the terms on the right in the theorem (i.e. in
$d'+d''$) are accounted for precisely once in the above enumeration of cases. 
\end{proof}

\ifx\undefined\bysame
\newcommand{\bysame}{\leavevmode\hbox to3em{\hrulefill}\,}
\fi

\end{document}